\documentclass[12pt,leqno]{amsart}

\numberwithin{equation}{section}

\newenvironment{proclaim}[1]
  {\par\medskip\noindent{\bf #1}  \it}
  {\rm\par\medskip}

\renewcommand{\a}{\alpha}

\newcommand{\C}{\Gamma}
\newcommand{\A}{\mathcal A}
\newcommand{\B}{\mathcal B}

\theoremstyle{plain}
  \newtheorem{thm}{Theorem}[section]
  \newtheorem{lem}[thm]{Lemma}
  \newtheorem{prop}[thm]{Proposition}
  \newtheorem{cor}[thm]{Corollary}

\theoremstyle{definition}
  \newtheorem{defn}[thm]{Definition}
  \newtheorem{example}[thm]{Example}

  \newtheorem{rem}[thm]{Remark}

\newtheorem*{example*}{Example}

\begin{document}

\title
[JSJ-decompositions and complexes of groups]
{JSJ-decompositions of finitely presented groups and
complexes of groups}

\author
[K. Fujiwara]{Koji Fujiwara}
\author
[P.Papasoglu]{Panos Papasoglu}

\dedicatory{Dedicated to Professor David Epstein
for his sixtieth birthday}


\email
[Koji Fujiwara]{fujiwara@math.tohoku.ac.jp}
\email
[Panos Papasoglu]{panos@math.uoa.gr}

\address
[Koji Fujiwara]
{Mathematical Institute, Tohoku University.
Sendai, 980-8578, Japan}
\address
[Panos Papasoglu] {Mathematics Department, University of Athens,
Athens 157 84, Greece }

\date{\today}

\keywords
{JSJ-decomposition, Complex of groups, Bass-Serre theory}
 \renewcommand{\subjclassname}{%
    \textup{2000} Mathematics Subject Classification}
\subjclass{20F65, 20E08, 57M07}

\begin{abstract}
  A JSJ-splitting of a group $G$ over a certain class of
subgroups is a graph of groups decomposition of $G$ which
describes all possible decompositions of $G$ as an amalgamated
product or an HNN extension over subgroups lying in the given
class. Such decompositions originated in 3-manifold topology. In
this paper we generalize the JSJ-splitting constructions of Sela,
Rips-Sela and Dunwoody-Sageev and we construct a JSJ-splitting for
any finitely presented group with respect to the class of all
slender subgroups along which the group splits. Our approach
relies on Haefliger's theory of group actions on CAT$(0)$ spaces.
\end{abstract}

\maketitle
\section{Introduction}
The type of graph of groups decompositions that we will consider
in this paper has its origin in 3-dimensional topology. Waldhausen
in \cite {W} defined the characteristic submanifold of a
3-manifold $M$ and used it in order to understand exotic homotopy
equivalences of 3-manifolds (i.e., homotopy equivalences that are
not homotopic to homeomorphisms). Here is a (weak) version of the
characteristic submanifold theory used by Waldhausen  that is of
interest to us: Let $M$ be a closed, irreducible, orientable
3-manifold. Then there is a finite collection of embedded 2-sided
incompressible tori such that each piece obtained by cutting $M$
along this collection of tori is either a Seifert fibered space or
atoroidal and acylindrical. Furthermore every embedded
incompressible torus of $M$ is  either homotopic to one of the
cutting tori or can be isotoped into a Seifert fibered space
piece. We note that embedded incompressible tori  of $M$
correspond to splittings of the fundamental group of $M$ over
abelian subgroups of rank 2. So from the algebraic point of view
we have a `description' of all splittings of $\pi _1(M)$ over
abelian groups of rank 2. Waldhausen in \cite {W} did not give a
proof of this theorem; it was proven later independently by
Jaco-Shalen \cite {JS}  and Johannson \cite {J} (this explains the
term JSJ-decomposition).

We recall that by Grushko's theorem every finitely generated group
$G$ can be decomposed as a free product of finitely many
indecomposable factors. Now if $G$ has no $\Bbb Z$ factors any
other free decomposition of $G$ is simply a product of a
rearrangement of conjugates of these indecomposable factors. One
can see JSJ-decomposition as a generalization of this description
for splittings of groups over certain classes of subgroups.

We recall that a group is termed {\it small} if it has no free subgroups
of rank 2.
Our paper deals with splittings over slender groups which are a
subclass of small groups. We recall that a finitely generated
group $G$ is {\it slender } if every subgroup of $G$ when it  acts
on a tree either leaves an infinite line invariant or it fixes a
point. It turns out that a group is slender if and only if all its
subgroups are finitely generated (see \cite {DS}). For example
finitely generated nilpotent groups are slender.

 To put our results on JSJ-decompositions
 in perspective we note that Dunwoody has shown that if $G$ is
a finitely presented group then if $\Gamma $ is a graph of groups
decomposition of $G$ with corresponding $G$-tree $T_{\Gamma }$
then there is a $G$-tree $T'$ and a $G$-equivariant map $\alpha
:T'\to T_{\Gamma }$ such that $T'/G$ has at most $\delta (G)$
essential vertices (see \cite{BF}, lemma 1). We recall that a
vertex in a graph of groups is not essential if it is adjacent to
exactly two edges and both edges and the vertex are labelled by
the same group. In other words one can obtain all graph of groups
decompositions of $G$ by `folding' from some graph of group
decompositions which have less than $\delta (G)$ vertices.

We remark that in general there is no bound on the number of
vertices of the graph of groups decompositions that one obtains
after folding. However in the special case of decompositions with
small edge groups Bestvina and Feighn (\cite {BF}. 
See Thm \ref{thm,bestvinafeighn} in this paper) have
strengthened this result showing that every reduced decomposition
$\Gamma $ of a finitely presented group $G$ with small edge groups
has at most $\gamma(G)$ vertices. Essentially they showed that in the
case of small splittings the number of `foldings' that keep the
edge groups small is bounded.
 The JSJ decomposition that we present here complements the
previous results as it gives a description of a set of
decompositions with slender edge groups from which we can obtain
any other decomposition by `foldings'. Roughly this set is
obtained as follows: we start with the JSJ decomposition and then
we refine it by picking for each enclosing group some splittings
that correspond to disjoint simple closed curves on the underlying
surface. Of course there are infinitely many such possible
refinements but they are completely described by the `surfaces'
that correspond to the enclosing groups.

Sela in \cite {S} was the first to introduce the notion of a
JSJ-decomposition for a generic class of groups, namely for
hyperbolic groups. Sela's JSJ-decomposition of hyperbolic groups
describes all splittings of a hyperbolic group over infinite
cyclic subgroups and was used to study the group of automorphisms
of a hyperbolic group. Sela's result was subsequently generalized
by Rips and Sela (\cite {RS}) to all finitely presented groups.
Dunwoody and Sageev (\cite {DS}) generalized this result further
and produced a JSJ-decomposition which describes all splittings of
a finitely presented group over slender groups under the
assumption that the group does not split over groups `smaller'
than the ones considered. Bowditch in \cite {B} gives a different
way of constructing the JSJ-decomposition of a hyperbolic group
using the boundary of the group. In particular this shows that the
JSJ-decomposition is invariant under quasi-isometries.

In this paper we produce for every finitely presented group $G$ a
JSJ-decomposition  of $G$ that describes all splittings of $G$
over all its slender subgroups.

Our approach to JSJ-decompositions differs from that of \cite
{S},\cite {RS} and of \cite {DS} in that we use neither $\Bbb
R$-trees nor presentation complexes. We use instead Haefliger's
theory of complexes of groups and actions on products of trees. To
see how this can be useful in studying splittings of groups
consider the following simple example: Let $G$ be the free abelian
group on two generators $a,b$. Then $G$ splits as an HNN extension
over infinitely many of its cyclic subgroups. Consider now two HNN
decompositions of $G$, namely the HNN decomposition of $G$ over
$\langle a \rangle$ and over $\langle b \rangle $. The trees corresponding to these
decompositions are infinite linear trees. Consider now the
diagonal action of $G$ on the product of these two trees. The
quotient is a torus. Every splitting of $G$ is now represented in
this quotient by a simple closed curve. We see therefore how we
can arrive at a description of infinitely many splittings by
considering an action on a product of trees corresponding to two
splittings.

Before stating our results we give a brief description of our
terminology: Let $T_A,T_B$ be Bass-Serre trees for one edge
splittings of a group $G$ over subgroups $A,B$. We say that the
splitting over $A$ is elliptic with respect to the splitting over
$B$ if $A$ fixes a vertex of $T_B$. If the splitting over $A$ is
not elliptic with respect to the splitting over $B$ we say that it
is hyperbolic. We say that the pair of two splittings is
hyperbolic-hyperbolic if they are hyperbolic with respect to each
other. We define similarly elliptic-elliptic etc (see
def.\ref{defn,type}). If a splitting over a slender group $A$ is
not hyperbolic-elliptic with respect to any other splitting over a
slender group then we say it is minimal. Finally we use the term
enclosing group (def.\ref{defn,enclosing}) for what Rips-Sela call
quadratically hanging group and Dunwoody-Sageev call hanging
${\mathcal K}$-by-orbifold group.

This paper is organized as follows. In section 2 we prove some
preliminary results and recall basic definitions from \cite {RS}.
In section 3 we introduce the notion of `minimality' of splittings
and prove several technical lemmas about minimal splittings that
are used in the sequel. In section 4 we apply Haefliger's theory
to produce `enclosing groups' for pairs of hyperbolic-hyperbolic
minimal splittings. Proposition \ref{prop,enclosinggroup} is the
main step in our construction of JSJ decompositions. It says that
we can always find a graph decomposition that contains both
splittings of a given pair of splittings. We note that, although
in our main theorem we consider only finitely presented groups,
proposition \ref{prop,enclosinggroup} is valid for groups that are
only finitely generated. Moreover proposition
\ref{prop,enclosinggroup} holds also for pairs of
hyperbolic-hyperbolic splittings over small groups.

 In section 5 using the
same machinery as in section 4 we show that there is a graph of
groups that `contains' all splittings from a family of
hyperbolic-hyperbolic minimal splittings (proposition \ref{prop
,enclosehyperbolic}). Using this we describe a refinement process
that produces the JSJ-decomposition of a finitely presented group
over all its slender subgroups. Because of the accessibility
results of Bestvina-Feighn (\cite{BF}. See Thm \ref{thm,bestvinafeighn}), 
there is an upper bound on
the complexity of graph decompositions that appear in the
refinement process, therefore this process must terminate. The
terminal graph decomposition must ``contain'' all minimal
splittings.

 The graph decomposition has special vertex groups (maybe
none) which are called maximal enclosing groups with adjacent edge
groups to be peripheral (see Def \ref{defn,enclosing}). Each of
them is an extension of the orbifold fundamental group of some
compact $2$-orbifold with boundary (maybe empty) by a slender
group, $F$. Examples are surface groups ($F$ is trivial) and the
fundamental group of a Seifert space ($F \simeq {\Bbb Z}$), which
is a $3$-manifold. We produce a graph decomposition using minimal
splittings (see Def \ref{defn,minimal}) of $G$. See Def
\ref{defn,type} for the definition of the type of a pair of
splittings, namely, hyperbolic-hyperbolic, elliptic-elliptic.

\begin{proclaim}{Theorem \ref{thm,jsj.minimal}.}
Let $G$ be a finitely presented group.
Then there exists a graph decomposition, $\Gamma$, of $G$
such that
\begin{enumerate}
\item
all edge groups are slender.
\item
Each edge of $\Gamma$ gives a
minimal splitting of $G$ along a slender group.
This splitting is elliptic-elliptic with respect
to any minimal splitting of $G$ along a slender subgroup.
\item
Each maximal enclosing group of $G$ is a conjugate
of some vertex group of $\Gamma$,
which we call a (maximal) enclosing vertex group.
The edge group of an edge adjacent to the vertex
of a maximal enclosing vertex group is a peripheral subgroup
of the enclosing group.

\item
Let $G=A*_C B$ or $A*_C$ be a minimal splitting along
a slender group $C$, and $T_C$ its Bass-Serre tree.
\begin{enumerate}
\item
If it is elliptic-elliptic with respect to all
minimal splittings of $G$ along slender groups, then
all vertex groups of $\Gamma$ are elliptic on $T_C$.
\item
If it is hyperbolic-hyperbolic
with respect to some minimal splitting of $G$ along a slender
group, then there is an enclosing vertex group, $S$, of $\Gamma$
which contains a conjugate of $C$,
which is unique among enclosing vertex groups of $\Gamma$.
$S$ is also the only one among
enclosing vertex groups which is hyperbolic on $T_C$.
There exist a base $2$-orbifold, $\Sigma$, for $S$ and
an essential simple closed curve or a segment on $\Sigma$
whose fundamental group (in the sense of
complex of groups) is a conjugate of $C$.

All vertex groups except for $S$ of $\Gamma$ are elliptic on $T_C$.

In particular, there is a graph decomposition, ${\mathcal S}$, of $S$
whose edge groups are in conjugates of $C$, which we can
substitute for $S$ in $\Gamma$ such that all vertex groups
of the resulting refinement of $\Gamma$ are elliptic on $T_C$.
\end{enumerate}
\end{enumerate}
\end{proclaim}

Although we produce $\Gamma$, called a JSJ-decomposition, using only minimal splittings,
it turns out that it is also good for non-minimal splittings.

\begin{proclaim}{Theorem \ref{thm,jsj}.}
Let $G$ be a finitely presented group, and $\Gamma$
 a graph decomposition we obtain in Theorem \ref{thm,jsj.minimal}.
Let $G=A*_C B, A*_C$ be a splitting along a slender group $C$, and
$T_C$ its Bass-Serre tree.
\begin{enumerate}
\item
If the group $C$ is elliptic with respect to any minimal splitting of $G$ along
a slender group, then all vertex groups of $\Gamma$ are
elliptic on $T_C$.
\item
Suppose the group $C$ is hyperbolic with respect to
some minimal splitting of $G$ along
a slender group. Then
\begin{enumerate}
\item all non-enclosing vertex groups of $\Gamma$ are elliptic on
$T_C$. \item For each enclosing vertex group, $V$, of $\Gamma$,
there is a graph decomposition of $V$, ${\mathcal V}$, whose edge
groups are in conjugates of $C$, which we can substitute for $V$
in $\Gamma$ such that if we substitute for all enclosing vertex
groups of $\Gamma$ then all vertex groups of the resulting
refinement of $\Gamma$ are elliptic on $T_C$.
\end{enumerate}
\end{enumerate}
\end{proclaim}
The first version of this paper is written in 1998. Since then a
very important application of JSJ-decompositions is found by
Z.Sela on Tarski's conjecture on the equivalence of the elementary
theory of $\mathbb{F}_2, \mathbb{F}_3$ (see \cite{S1} and the
following papers of Sela on this). He uses JSJ-decompositions
along abelian subgroups. We note also that the question of
`uniqueness' of JSJ-splittings has been treated in \cite{Fo}. We
would like to thank M.Bestvina, M.Feighn, V.Guirardel, B.Leeb, M.Sageev, Z.Sela
and G.A.Swarup for discussions related to this work. We would like
to thank A. Haefliger for his interest in this work and many
suggestions that improved the exposition. Finally we would like to
thank the referee for detailed suggestions which we found very
helpful.

\section{Pairs of splittings}
In this section we recall and generalize notation from \cite {RS}.
\begin{defn}[types of a pair]\label{defn,type}
 Let $A_1\star _{C_1}B_1$ (or $A_1\star _{C_1}$),
 $A_2\star _{C_2}B_2$ (or $A_2\star _{C_2}$) be two splittings
of a finitely generated  group $G$ with corresponding Bass-Serre
trees $T_1$, $T_2$.
We say that the first splitting is hyperbolic with respect to the second
if there is $c_1\in C_1$
acting as a hyperbolic element on $T_2$. We say that the first splitting
is elliptic
with respect to the second if $C_1$ fixes a point of $T_2$.
We say that this pair
of splittings is hyperbolic-hyperbolic if each splitting is hyperbolic
with respect to the other.
Similarly we define what it means for a pair of splittings to be
elliptic-elliptic, elliptic-hyperbolic
and hyperbolic-elliptic.
\end{defn}

It is often useful to keep in mind the `geometric' meaning of this
definition: Consider
for example a closed surface. Splittings of its fundamental group over
$\Bbb Z$
 correspond to simple closed curves on the surface. Two splittings are
hyperbolic-hyperbolic
if their corresponding curves intersect and elliptic-elliptic otherwise.
Consider now a punctured surface  and two splittings of its fundamental
group: one corresponding to
a simple closed curve (a splitting over $\Bbb Z$) and a free splitting
corresponding to
an arc having its endpoints on the puncture such that the two curves
intersect at one point. This
pair of splittings is hyperbolic-elliptic.

\begin{prop}
\label{prop,hyp-hyp}
Let $A_1\star _{C_1}B_1$ (or $A_1\star _{C_1}$), $A_2\star _{C_2}B_2$
(or  $A_2\star _{C_2}$) be two splittings
of a group $G$ with corresponding Bass-Serre trees $T_1$, $T_2$. Suppose
that
there is no splitting of $G$ of the form $A\star _{C}B$ or $A\star _{C}$
 with $C$ an infinite index subgroup of $C_1$ or of $C_2$.
Then this pair of splittings is either hyperbolic-hyperbolic or
elliptic-elliptic.
\end{prop}

\proof We treat first the amalgamated product case. Let $T_1$,
$T_2$ be the Bass-Serre trees of the two splittings $A_1\star
_{C_1}B_1$, $A_2\star _{C_2}B_2$. Suppose that $C_1$ does not fix
any vertex of $T_2$ and that $C_2$ does fix a vertex of $T_1$.
Without loss of generality we can assume that $C_2$ fixes the
vertex stabilized by $A_1$. Consider the actions of $A_2,B_2$ on
$T_1$. Suppose that both $A_2,B_2$ fix a vertex. If they fix
different vertices then $C_2$ fixes an edge, so it is a finite
index subgroup of a conjugate of $C_1$. But then $C_1$ can not be
hyperbolic with respect to $A_2\star _{C_2}B_2$.
 On the other hand it is not possible that they fix the same vertex
since $A_2,B_2$ generate $G$. So at least one of them, say $A_2$,
does not fix a vertex. But then the action of $A_2$ on $T_1$
induces a splitting of $A_2$ over a group $C$ which is an infinite
index subgroup of $C_1$. Since $C_2$ is contained in a vertex
group of this splitting we obtain a splitting of $G$ over $C$
which is a contradiction.

We consider now the case one of the splittings is an
HNN-extension: say  we have the splittings $A_1\star _{C_1}B_1$,
$A_2\star _{C_2}$ with Bass-Serre trees $T_1,T_2$. Assume $C_1$ is
hyperbolic on $T_2$ and $C_2$ elliptic on $T_1$. Again it is not
possible that $A_2$ fix a vertex of $T_1$.  Indeed $C_2=A_2\cap
tA_2t^{-1}$ and if $A_2$ fixes a vertex $C_2$ is contained in a
conjugate of $C_1$ which is impossible (note that $t$ can not fix
the same vertex as $A_2$). We can therefore obtain a splitting of
$G$ over an infinite index subgroup of $C_1$ which is a
contradiction. If $C_1$ is elliptic on $T_2$ and $C_2$ hyperbolic
on $T_1$ we argue as in the first case. The case where both
splitting are HNN extension is treated similarly. \qed

\begin{rem}\label{rem,infiniteindex}
 In the proof of the Proposition \ref{prop,hyp-hyp} one shows
in fact that
if $A_2\star _{C_2}B_2$ is elliptic with respect to
$A_1\star _{C_1}B_1$
then either
$A_1\star _{C_1}B_1$ is elliptic too, or there is a splitting of $G$
over a subgroup
of infinite index of $C_1$.
\end{rem}

\section{Minimal splittings}
\begin{defn}[Minimal splittings]\label{defn,minimal}
We call a splitting $A\star _{C}B$ (or $A\star _{C}$) of a group
$G$ minimal if it is not hyperbolic-elliptic with respect to any
other splitting of $G$ over a slender subgroup.
\end{defn}

\begin{rem}\label{rem,example}
Remark \ref{rem,infiniteindex}
 implies that if $G$ splits over $C$ but does not
split over an
infinite index subgroup of $C$ then the splitting of $G$ over $C$ is
minimal.
There are examples of minimal and non-minimal splittings
over a common subgroup.
For example let $H$ be a group which does not
split and let $G={\Bbb Z} ^2 *H$. If $a,b$ are generators of
${\Bbb Z}^2$ the splitting
of ${\Bbb Z}^2$ over $\langle a \rangle$ induces a minimal splitting of $G$ over
$\langle a \rangle $.
On the
other hand the splitting of $G$ given by $G={\Bbb Z}^2 * _{\langle a \rangle}
(H* \langle a \rangle )$
is not
minimal. Indeed it is hyperbolic-elliptic with respect to the splitting
of $G$
over $\langle b \rangle $ which is induced from the splitting of  $\Bbb Z ^2$ over
$\langle b \rangle $.
\end{rem}

We collect results on minimal splittings we need.
We first show the following:

\begin{lem}
\label{lem,preminimal}
Suppose that
a group G splits over the slender groups $C_1,C_2$ and $K\subset C_2$.
Assume moreover
that the splittings  over $C_1,C_2$ are
hyperbolic-hyperbolic,
the splitting over $C_1$ is minimal and that $G$
admits an action on a
tree $T$ such that $C_2$ acts hyperbolically and $K$ fixes a vertex.
Then the splittings over $C_1$ and $K$ are not
hyperbolic-hyperbolic.
\end{lem}

\proof
We will prove this by contradiction.
Let $T_1,T_2,T_3$ be, respectively, the Bass-Serre trees of
the splittings over $C_1,C_2,K$.
Without loss of generality
we can assume the axes of $C_2,K$ when acting on $T_1$ contain
an edge stabilized by $C_1 $.
Let $t\in K\subset C_2$ be an element acting hyperbolically
on $T_1$. Similarly let $u\in C_1$ be
an element acting hyperbolically on $T_2,T_3$ and  $y\in C_2$ be an
element
acting hyperbolically on $T$.
We distinguish 2 cases:

\noindent
{\it Case 1}: $y$ acts elliptically on $T_1$.  Then either $y$ fixes
the axis of translation of $C_2$ or it acts on it by a reflection (in
the
dihedral action case). In both cases $y^2\in C_1$.
Since $y ^2\notin K$ we have that $y^2$ acts hyperbolically on $T_3$.
Therefore
there are $m,n\in \Bbb Z$ such that $y^mu^n $ fixes an edge of $T_3$ and
$y^mu^n\in C_1\cap zKz^{-1} $, hence
$y^mu^n \in C_1\cap xC_2x^{-1}  $. This is clearly impossible
since $y$ is elliptic
when acting on $T_2$ while $u$ is hyperbolic.
\newline
{\it Case 2}: $y$ acts hyperbolically on $T_1$.  Without loss of
generality we assume that $t$ fixes the axis of translation of $y$ on
$T$.
Indeed if this is not so we can replace $t$ by $t^2$.
Since both $t$ and $y$ act hyperbolically on $T_1$
 there are
$m,n\in \Bbb Z$ such that
$t^my^n\in C_1$. On the other hand $t^my^n$ acts hyperbolically on $T$
since $t$ fixes the axis of translation of $y$.
So $t^my^n$ does not lie in a conjugate of $K$.
For the same reason $(t^my^n)^2$ does not lie in a conjugate of $K$.
We consider now the action of $C_1$ on $T_3$. If  $t^my^n$ is elliptic
then $(t^my^n)^2$ fix  the axis of translation of $C_1$ and therefore
lies in a conjugate of $K$, which is impossible as we noted above.
Therefore
both  $t^my^n$  and $u$ acts hyperbolically on $T_3$.
We infer that there are $p,q\in \Bbb Z$ such that
$ (t^my^n)^pu^q$ lies in a conjugate of $K$. Therefore this element
fixes
the translation axis of $C_1$ when acting on $T_2$. This is however
impossible since $t^my^n\in C_1\cap C_2$ so it fixes the axis while
$u$ acts hyperbolically on $T_2$.
This finishes the proof of the lemma.
\qed

Using lemma \ref{lem,preminimal}, we show the following.

\begin{prop}[dual-minimality]
\label{prop,dual-minimal}
Let $A_1\star _{C_1}B_1$ (or $A_1\star _{C_1}$) be a minimal
splitting of $G$ over a slender group $C_1$. Suppose that
$A_1\star _{C_1}B_1$ ( $A_1\star _{C_1}$)is hyperbolic-hyperbolic
with respect to another splitting of $G$, $A_2\star _{C_2}B_2$ (or
$A_2\star _{C_2}$), where $C_2$ is slender. Then $A_2\star
_{C_2}B_2$ (or $A_2\star _{C_2}$) is also minimal.
\end{prop}

\proof
We denote the Bass-Serre trees  for the splittings over $C_1,C_2$
by $T_1,T_2$ respectively.

 Suppose that
the splitting over $C_2$ is not minimal; then it is
hyperbolic-elliptic with respect to another splitting over a
slender subgroups $C_3$. We distinguish 2 cases:

\textit{1st case}: The splitting over $C_3$ is an amalgamated
product, say $A_3\star _{C_3}B_3$.

 We let
$A_3,B_3$ act on $T_2$ and we get graph of groups decompositions
for $A_3,B_3$, say $\Gamma _1,\Gamma _2$. Since $C_3$ is elliptic
when acting on $T_2$ we can refine $A_3\star _{C_3}B_3$ by replacing
$A_3,B_3$ by $\Gamma _1,\Gamma _2$. We collapse then the
edge labelled by $C_3$ and we obtain a new graph of groups decomposition
that we call $\Gamma $. We note that all vertex groups of $\Gamma $ fix
vertices of $T_2$. This implies that $C_1$ is not contained in
a conjugate of a vertex group of $\Gamma $. Therefore
we can collapse all edges of $\Gamma $ except one and obtain
a splitting over a subgroup $K$ of $C_2$ such that $C_1$ is hyperbolic
with respect to this splitting. Since $C_1$ is minimal the pair of
splittings
over $C_1,K$ is
hyperbolic-hyperbolic.
 This however contradicts lemma \ref{lem,preminimal}
since $K$ fixes a vertex of the Bass-Serre tree of $A_3\star _{C_3}B_3$
while $C_2$ acts hyperbolically on this tree.

\textit{2nd case}: The splitting over $C_3$ is an HNN-extension,
say $A_3\star _{C_3}$.

The argument is similar in this case but a bit more delicate. We
let $A_3$ act on $T_2$ and we obtain a graph decomposition for
$A_3$, say $\Gamma _1$. Since $C_3$ fixes a vertex of $T_2$ we can
refine $A_3\star _{C_3}$ by replacing $A_3$ by $\Gamma _1$. Let
$e$ be the edge of $A_3\star _{C_3}$. If $e$ stays a loop in
$\Gamma _1$ we argue as in the amalgamated product case. After we
collapse $e$ in $\Gamma _1$ we obtain a graph of groups such that
$C_1$ is not contained in the conjugate of any vertex group. We
arrive then at a contradiction as before.

Assume now that $e$ connects to two different vertices  in $\Gamma
_1$. Let $V,U$ be the vertex groups of these vertices. Clearly
$C_2$ is not contained in a conjugate of either $V$ or $U$. We
remark that the vertex we get after collapsing $e$ in $\Gamma _1$
is labelled by $\langle V_1,V_2 \rangle $. By Bass-Serre theory $C_2$ is not
contained in a conjugate of $\langle V_1,V_2 \rangle $ either. Let's call $\Gamma
$ the graph of groups obtained after collapsing $e$. Let
$T_{\Gamma } $ be its Bass-Serre tree. Clearly $C_2$ does not fix
any vertex of $T_{\Gamma } $. It follows that we can collapse all
edges of $\Gamma $ except one and obtain an elliptic-hyperbolic
splitting with respect to $C_2$. What we have gained is that this
splitting is over a subgroup of $C_2$, say $K$. Of course if this
splitting is an amalgam we are done by case 1 so we assume it is
an HNN-extension $A\star _K$. Let's call $e_1$ the edge of this
HNN-extension and let $T_K$ be its Bass-Serre tree. If $C_1$ acts
hyperbolically on $T$ then we are done as before by lemma
\ref{lem,preminimal} since $C_2$ is hyperbolic on $T$ and $K$
elliptic. Otherwise we let $A$ act on $T_2$ and we refine $A\star
_K$ as before. Let's call $\Gamma '$ the graph of groups obtained.
If $e_1$ stays a loop in $\Gamma '$ we are done as before.
Otherwise by collapsing $e_1$ we obtain a graph of groups such
that no vertex group contains a conjugate of  $C_1$. We can
collapse this graph further to a one edge splitting over, say
$K_1<C_2$ such that its vertex groups do not contain a conjugate
of $C_1$. Since the splitting over $C_1$ is minimal this new
splitting is hyperbolic-hyperbolic with respect to the splitting
over $C_1$. Moreover $K_1$ fixes a vertex of $T$ while $C_2$ acts
hyperbolically on $T$. This contradicts lemma
\ref{lem,preminimal}.
\qed

We prove an accessibility result for minimal splittings.

\begin{prop}[accessibility of minimal splittings]
\label{prop,accessibility}
 Let $G$ be a finitely generated group. There is no
infinite sequence of splittings of $G$ of the form $A_n\star _{C_n}B_n$
or of the form  $A_n\star _{C_n}$
where $C_{n+1}$ is a subgroup of $C_n$, $C_1$ is a (finitely generated)
slender group, such that
 $C_n$ acts hyperbolically on some $G$-tree $T_n$ while
$C_{n+1}$ fixes a vertex of $T_n$.
\end{prop}

\proof
We define a sequence of homomorphisms $f_n$ from $C_1$ to ${\Bbb Z}/2{\Bbb Z}$
as follows:
Consider the graph of groups corresponding to the action of $C_1$ on $T_n$.
The fundamental
group of this graph of groups is $C_1$. If the underlying graph of
this graph of groups is a circle,
map this group to ${\Bbb Z}/2{\Bbb Z}$ by mapping all vertex
groups to $0$ and the single loop of the graph to $1$.
If $C_1$ acts by a dihedral  type action on $T_n$,
map $C_1$ to ${\Bbb Z}/2{\Bbb Z} *{\Bbb Z} /2{\Bbb Z}$ in the obvious way and then
map ${\Bbb Z}/2{\Bbb Z} *{\Bbb Z} /2{\Bbb Z}$ to ${\Bbb Z} /2{\Bbb Z}$
so that $C_{n+1}$
is mapped to $0$. This is possible since $C_{n+1}$ is elliptic.
By construction
$f_{n+1}(C_{n+1})={\Bbb Z} /2{\Bbb Z}$. It follows
that the map $\Phi_n:C_1\to ({\Bbb Z} /2{\Bbb Z})^n$
is onto for every $n$. This contradicts the fact that $C_1$ is finitely
generated.
\qed

Each edge, $e$, of a graph decomposition, $\Gamma$, of a group
$G$ gives (rise to) a splitting of $G$ along the edge group
of $e$, $E$, by collapsing all edges of $\Gamma$
but $e$. We do this often in the paper.
To state the main result of this section, we give
one definition.

\begin{defn}[Refinement]
\label{defn,refinement}
Let $\Gamma $ be a graph of groups decomposition of $G$. We say that
$\Gamma '$ is a refinement of $\Gamma $ if each vertex group of $\Gamma '$ is
contained in a conjugate of a vertex group of $\Gamma $. We say that $\Gamma '$
is a proper
refinement of $\Gamma $ if $\Gamma '$ is a refinement of $\Gamma $ and $\Gamma $
is not a refinement of $\Gamma '$.
\end{defn}

Let $\Gamma $ be a graph of groups and let $V$ be a vertex group
of $\Gamma $, which admits a graph of groups decomposition $\Delta
$ such that a conjugate of each edge group adjacent to $V$ in
$\Gamma$ is contained in a vertex group of $\Delta $. Then one can
obtain a refinement of $\Gamma $ by replacing $V$ by $\Delta $. In
the refinement, each edge, $e$, in $\Gamma$ adjacent to the vertex
for $V$ is connected to a vertex of $\Delta$ whose vertex group
contains a conjugate of the edge group of $e$. The monomorphism
from the edge group of $e$ to the vertex group of $\Delta $ is
equal to the corresponding monomorphism in $\Gamma $ modified by
conjugation.

This is a special type of refinement used often in this paper. We
say that we {\it substitute} $\Delta$ for $V$ in $\Gamma$.

\begin{prop}[Modification to minimal splittings]
\label{prop,modification} Let $G$ be a finitely presented group.
Suppose that $\Gamma $ is a graph of groups decomposition of $G$
with slender edge groups. Then there is a graph of groups
decomposition of $G$, $\Gamma '$, which is a refinement of $\Gamma
$ such that all edges of $\Gamma '$ give rise to minimal
splittings of $G$. All edge groups of $\Gamma'$ are subgroups of
edge groups of $\Gamma$.
\end{prop}

\proof
We give two proofs of this proposition. We think that in the first
one the idea is more transparent, which uses actions on product of
trees. Since we use terminologies and ideas from the part we
construct "enclosing groups" in Proposition
\ref{prop,enclosinggroup}, one should read the first proof after
reading that part. The second proof uses only classical Bass-Serre
theory and might be more palatable to readers not accustomed to
Haefliger's theory.
\newline
{\it 1st proof}.
We define a process to produce a
sequence
of refinements of $\C$ which we can continue
as long as an edge of a graph decomposition in the sequence gives a non-minimal splitting,
and then show  that
it must terminate in a finite step.

If every edge of $\C$ gives a minimal splitting, nothing to do, so suppose there is
an edge, $e$, of $\C$ with the edge group $E$ which gives a splitting
which is is hyperbolic-elliptic with respect to, say, $G=P*_R Q$ or $P*_R$.
We consider the action of $G$ on the product of trees
$T_{\C},T_R$ where $T_{\C}$ is the Bass-Serre tree of $\C $ and $T_R$ the tree
of the splitting over $R$.
We consider the diagonal action of $G$ on $T_{\C}\times T_R$,
then produce a $G$-invariant subcomplex of $T_{\C}\times T_R$, $Y$,
such that $Y/G$ is compact as we will do in the proof of
Proposition \ref{prop,enclosinggroup}.
$Y/G$ has a structure of a complex of groups whose fundamental
group is $G$.
In the construction of $Y/G$, we give priority (see Remark \ref{rem,priority})
to the decomposition $\C$ over the splitting of $G$ along $R$, so that we can
recover $\C$ from $Y/G$ by collapsing the complex of groups obtained. To fix ideas
we think of $T_{\C}$ as horizontal (see the paragraph after
Lemma \ref{lem,vankampen}).
Since the slender group $E$ acts hyperbolically
on $T_R$, there is a line, $l_E$, in $T_R$ which
is invariant by $E$. Then, $l_E/E$ is either
a segment, when the action of $E$ is dihedral, or else, a circle.
Put $c_E=l_E/E$, which is the core for $E$.
There is a map from $c_E \times [0,1]$ to $Y/G$, and
let's call the image, $b_E$, the band for $E$.
$b_E$ is a union of finite squares.
For other edges, $e_i$, of $\Gamma$ than $e$ with
edge groups $E_i$, we have similar objects, $b_{E_i}$,
which can be a segment, when the action of $E_i$ on $T_R$
is elliptic. As in Proposition \ref{prop,enclosinggroup},
the squares in $Y/G$ is exactly the union of the squares
contained in $b_E$ and other $b_{E_i}$'s.
Note that there is at least one square in $Y/G$, which
is contained in $b_E$.

On the other hand since $R$ is elliptic on $T_{\Gamma}$,
$T_{\Gamma}/R$ is a tree, so that the intersection
of $T_{\Gamma}/R$ (note that this is naturally embedded in
$T_{\C}\times T_R$) and $Y/G$ is a forest.
Therefore, we can remove squares from $Y/G$ without
changing the fundamental group, which is $G$.
We then collapse all vertical edges which are left.
In this way, we obtain a graph decomposition of
$G$, which we denote $\C _1$.
By construction, all vertex groups of $\Gamma_1$ are
elliptic on $T_{\Gamma}$, so that
$\Gamma_1$ is a refinement of $\Gamma$.
Also edge groups of $\Gamma_1$ are subgroups of
edge groups of $\Gamma$.
There is no edge group of $\C _1$ which is hyperbolic-elliptic with
respect to the splitting over $R$.
If all edge of $\Gamma_1$ gives a minimal splitting of $G$, then
$\Gamma_1$ is a desired one. If not, we apply the same
process to $\Gamma_1$, and obtain a refinement, $\Gamma_2$.
But this process must terminate by Prop
\ref{prop,accessibility}
and Theorem \ref{thm,bestvinafeighn}, which gives a desired one.
\newline
{\it 2nd proof}.
If all edges of $\Gamma $ correspond to minimal splittings then there is nothing to prove.
Assume therefore that an edge $e$ of $\Gamma $ corresponds to a splitting which is not
minimal. We will construct a finite sequence of refinements of $\Gamma $ such that
the last term of the sequence is $\Gamma '$.
 Let's say that $e$ is labelled by the slender group $E$. Let
$A*_EB$ (or $A*_E$) be the decomposition of $G$ obtained by collapsing all edges
of $\Gamma $ except $e$. Since this splitting is not minimal it is hyperbolic-elliptic
with respect to another
splitting of $G$ over a slender group, say $P*_RQ$ (or $P*_R$).

Let $T_E,T_R$
 be the Bass-Serre tree of the splittings
 of $G$ over $E,R$
 and let $T_{\Gamma }$ be
the Bass-Serre tree corresponding to $\Gamma $.
We distinguish two cases:\newline
Case 1: $R$ is contained in a conjugate of $E$.\newline
In this case we let $P,Q$ act on $T_{\Gamma }$. We obtain graph of groups
decompositions of $P,Q$ and we refine $P*_RQ$ by substituting $P,Q$ by these
graphs of groups decompositions. In this way we obtain a graph of groups
decomposition $\Gamma _1$. If some edges of $\Gamma _1$ (which are not loops) are
labelled by the same group as an adjacent vertex we collapse them. For simplicity
we still call $\Gamma _1$ the graph of groups obtained after this collapsing.
We remark that all edge groups
of $\Gamma _1$ fix a vertex of the tree of the splitting $P*_RQ$. We argue in the same way
if the spitting over $R$ is an HNN-extension.\newline
Case 2: $R$ is not contained in any conjugate of $E$. \newline
We let $P,Q$ act on $T_E$ and we obtain graph of groups decompositions of these groups.
We refine $P*_RQ$ as before by substituting $P,Q$ by the graph of groups obtained.
We note that by our hypothesis in case 2 $P,Q$ do not fix both vertices of $T_E$.
Let's call $\Delta $ the graph of groups obtained in this way. Since $C$ fixes a vertex
of $T_E$ we can assume without loss of generality that $C\subset P$.
 We collapse
the edge of $\Delta $ labelled by $R$ and we obtain a graph of groups $\Delta _1$.
We note now that if $E$ acts hyperbolically on the Bass-Serre tree of $\Delta _1$
we are in the case 1 (i.e., we have a pair of
hyperbolic-elliptic splittings where the second splitting is obtained by appropriately
collapsing all edges of $\Delta _1$ except one). So we can refine $\Gamma $ and
obtain a decomposition $\Gamma _1$ as in case 1. We suppose now
that this is not the case.
Let's denote by $P'$ the group of the vertex obtained after collapsing the edge labelled by $R$.
We let $P'$ act on $T_{\Gamma }$ and we obtain a graph of groups decomposition
of $P'$, say $\Delta _2$.
We note that the vertex obtained after this collapsing is now labelled
by a subgroup of $P$, say $P'$.
We let $P'$ act on $T_E$ and we obtain a graph of groups decomposition of $P'$, say $\Delta _2$.
If every edge group of $\Delta _1$ acts elliptically on $T_R$ then we let all other
vertices of $\Delta _1$ act on $T_{\Gamma }$ and we substitute all these vertices
in $\Delta _1$ by the graphs obtained. We also substitute $P'$ by $\Delta _2$.
We call the graph of groups obtained in this way $\Gamma _1$.

Finally we explain what we do if some edge of $\Delta _2$ acts hyperbolically on $T_R$.
We note that $P'$ splits over $R$, indeed $P'$ corresponds to a one-edge subgraph of $\Delta $.
Abusing notation we call still $T_R$ the tree of the splitting of $P'$ over $R$.
We now repeat with $P'$ the procedure applied to $G$. We note that we are necessarily
in case 2 as $R$ can not be contained in a conjugate of an edge group of $\Delta _2$.
As before we either obtain a refinement of $\Delta _2$ such that all edge groups of $\Delta _2$
act elliptically on $T_R$ or we obtain a non-trivial decomposition of $P'$, say $\Delta '$,
such that an edge of $\Delta '$ is labelled by
$R$ and the following holds: If we collapse
the edge of $\Delta '$ labelled by $R$ we obtain a vertex $P"$ which has the same property
as $E'$. Namely if $\Delta _3$ is the decomposition of $P"$ obtained by acting on $T_{\Gamma }$
then some edge group of $\Delta _3$ acts hyperbolically on $T_R$. If we denote by $\Delta '_1$
the decomposition of $P'$ obtained after the collapsing we remark that we can substitute
$P'$ by $\Delta '_1$ in $\Delta _1$ and obtain a decomposition of $G$
with more edges than $\Delta _1$.
Now we repeat the same procedure to $P"$.
By Theorem \ref{thm,bestvinafeighn},
this process terminates and produces a refinement of $\Gamma $ which we call $\Gamma _1$.

By the argument above we obtain in both cases a graph of groups decomposition of $G$ $\Gamma _1$
which has the following properties:\newline
1) $\Gamma _1$ is a refinement of $\Gamma $ and \newline
2) There is an action of $G$ on a tree $T$ such that some edge group of $\Gamma $ act
on $T$ hyperbolically while all edge groups of $\Gamma _1$ act on $T$ elliptically.

Now we repeat the same procedure to $\Gamma _1$ and we obtain a graph of groups $\Gamma _2$ etc.
One sees that this
 procedure will terminate
using Prop \ref{prop,accessibility}
and Theorem \ref{thm,bestvinafeighn}.
The last step of this procedure
 produces
 a decomposition $\Gamma '$ as required by this proposition.
\qed

One finds an argument similar to the 1st proof, using product of trees and
retraction in the paper \cite{DF}.

\section{Enclosing groups for a pair of hyperbolic-hyperbolic splittings}
\subsection{Product of trees and core}\label{ss,haefliger}
Producing a graph of groups which `contains' a given pair of
hyperbolic-hyperbolic splittings along slender groups is the main
step in the construction of a JSJ-decomposition of a group. This
step explains also what type of groups should appear as vertex
groups in a JSJ-decomposition of a group. We produce such a graph
of groups in proposition \ref{prop,enclosinggroup}.

 We
recall here the definitions of a complex of groups and the
fundamental group of such a complex. They were first given in the
case of simplicial complexes in \cite {H} and then generalized to
polyhedral complexes in \cite {BH}. Here we will give the
definition only in the case of 2-dimensional complexes. We
recommend \cite {BH} ch.III.C for a more extensive treatment.

 Let $X$ be a polyhedral complex of dimension less or equal to 2.

 We associate to $X$ an oriented graph as follows: The vertex set
 $V(X)$ is the set of $n$-cells of $X$ (where $n=0,1,2$).
 The set of oriented edges $E(X)$ is the set
 $E(X)=\{(\tau, \sigma) \}$ where $\sigma $ is an $n$-cell of
 $X$ and $\tau $ is a face of $\sigma $. If $e\in E(X)$, $e=(\tau,
 \sigma)$, we define the original vertex of $e$, $i(e)$ to be $\tau $ and
 the terminal vertex of $e$, $t(e)$ to be $\sigma $.

  If $a,b\in
 E(X)$ are such that $i(a)=t(b)$ we define the composition $ab$ of $a,b$ to
 be the edge $ab=(i(b),t(a))$. If $t(b)=i(a)$ for edges $a,b$ we say
 that these edges are composable.
 We remark that the set $E(X)$ is in fact the set of edges of the
 barycentric subdivision of $X$. Geometrically one represents the
 edge $e=(\tau,\sigma)$ by an edge joining the barycenter of $\tau
 $ to the barycenter of $\sigma $. Also $V(X)$ can be identified with the set of vertices
of the barycentric subdivision of $X$, to a cell $\sigma $ there
corresponds a vertex of the barycentric subdivision, the
barycenter of $\sigma $.

 A complex of groups $G(X)=(X,G_{\sigma },\psi _{a },g_{a,b})$
 with underlying complex $X$ is given by the following data: \newline
 1. For each
 $n$-cell of $X$, $\sigma $ we are given a group
 $G_{\sigma}$.\newline
 2. If $a$ is an edge in $E(X)$ with $i(a)=\sigma $, $t(a)=\tau $
 we are given an injective homomorphism $\psi _a:G_{\sigma }\to
 G_{\tau }$. \newline
3. If $a,b$ are composable edges we are given an element
$g_{a,b}\in E_{t(a)}$ such that $$g_{a,b}\psi
_{ab}g_{a,b}^{-1}=\psi_a\psi_b $$

We remark that when $\dim (X)=1$, $G(X)$ is simply a graph of
groups. In fact in Haefliger's setup loops are not allowed so to
represent a graph of groups with underlying graph $\Gamma $ one
eliminates loops by passing to the barycentric subdivision of
$\Gamma $. In this case there are no composable edges so condition
3 is void.

We define the fundamental group of a complex of groups
$\pi_1(G(X),\sigma _0)$ as follows:

Let $E^{\pm}(X)$ be the set of symbols $a^+$, $a^-$ where $a\in
E(X)$. Let $T$ be a maximal tree of the graph $(V(X),E(X))$
defined above. $\pi_1(G(X),\sigma _0)$ is the group with
generating set:
$$\coprod G_{\sigma },\, \sigma \in V(X), \, \coprod E^{\pm}(X) $$
and set of relations: \smallskip

$ relations \ \  of \ \ G_{\sigma },\, (a^+)^{-1}=a^-,\,
(a^-)^{-1}=a^+,\, (\forall a \in E(X))$ \smallskip

$ a^+b^+=g_{a,b}(ab)^+, \forall a,b \in E(X),\, \psi
_a(g)=a^+ga^-,\, \forall g\in G_{i(a)},\, a^+=1,\, \forall a\in T$
\bigskip

It is shown in \cite {BH} that this group does not depend up to
isomorphism on the choice of maximal tree $T$ and its elements can
be represented by `homotopy classes' of loops in a similar way as
for graphs of groups.

It will be useful for us to define barycentric subdivisions of
complexes of groups $G(X)$. This will be an operation that leaves
the fundamental group of the complex of groups unchanged but
substitutes the underlying complex $X$ with its barycentric
subdivision $X'$. We explain this first in the case of graphs of
groups. If $G(X)$ is a graph of groups then we have a group $G_v$
for each vertex $v$ of the barycentric subdivision of $X$. If $v$
is the barycenter of an edge $e$, $G_v$ is by definition in
Haefliger's notation the group associated to the edge $e$, $G_e$.
Now if $v$ is a vertex of the second barycentric subdivision then
$v$ lies in some edge $e$ of $X$ so we define $G_v=G_e$. The
oriented edges $E(X)$ are of two types:

1) an edge $e$ from a barycenter $v$ of the second barycentric
subdivision to a barycenter $w$ of the first barycentric
subdivision. If this case the map $\psi _e: G_v\to G_w$ is the
identity.

2) an edge $e$ from a barycenter $v$ of the second barycentric
subdivision to a vertex $w$ of $X$. In this case $v$ is the
barycenter of an edge $a$ of the first barycentric subdivision and
$G_v$ is isomorphic to $G_{i(a)}$. We define then $\psi _e:G_v\to
G_w$ to be $\psi _e$.

Let's call $G(X')$ the graph of groups obtained by this operation.
It is clear that the fundamental group of $G(X')$ is isomorphic to
the fundamental group of $G(X)$.

Let now $X$ be a 2-dimensional complex and $G(X)$ a complex of
groups with underlying complex $X$. Let $X'$ be the barycentric
subdivision of $X$. We associate to $X'$ a graph $((V(X'),E(X'))$
as we did for $X$. Now we explain what are the groups and maps
associated to $(V(X'),E(X'))$.

In order to describe the groups associated to $V(X')$ it is
convenient to recall the geometric representation of $V(X),E(X)$.

The vertices of $V(X)$ correspond geometrically to barycenters of
$n$-cells of $X$, i.e. they are just the vertices of the
barycentric subdivision of $X$. Similarly the edges $E(X)$ are the
edges of the barycentric subdivision and the orientation of an
edge is from the barycenter of a face of $X$ to a vertex of $X$.

$V(X')$ analogously can be identified with the set of vertices of
the second barycentric subdivision of $X$ and the edges $E(X')$
with the edge set of the second barycentric subdivision of $X$.
All 2-cells of the barycentric subdivision of $X$ are 2-simplices.

If $v$ is a vertex of $V(X')$ which is the barycenter of the
2-simplex $\sigma $ then there is a single 2-cell $\tau $ of $X$
containing $\sigma $. If $w$ is the barycenter of $\tau $ we
define $G_v=G_w$.
 If
$v$ is a vertex of $V(X')$ which is the barycenter of an edge $a$
we define $G_{v }=G_{i(a)}$.

We explain now what are the homomorphisms corresponding to
$E(X')$. If $a$ is an edge of $E(X')$ then there are two cases:

1) $i(a)$ is the barycenter of a 2-simplex $\sigma $ of $X'$. Then
if $t(a)$ is the barycenter of a 2-cell $\tau $ of $X$ by
definition $G_{i(a)}=G_{t(a)}$ and we define $\psi _a$ to be the
identity. Otherwise if $w$ is the barycenter of the 2-cell of $X$
containing $\sigma $ we have that $G_{i(a)}=G_{w}$ and there is an
edge $e$ in $E(X)$ from $w$ to $t(a)$. We define then $\psi _a$ to
be $\psi _e$.

2) $i(a)$ is the barycenter of an edge $e$ of $X'$. If $t(a)=i(e)$
we define $\psi _a$ to be the identity. Otherwise $t(a)=t(e)$ and
we define $\psi _a$ to be $\psi _e$.

It remains to define the `twisting elements' for pairs of
composable edges of $G(X')$. We remark that if $a',b'$ are
composable edges of $X'$ then either $\psi _{b'}=id$ and $\psi
_{a'b'}=\psi _{a'}$ in which case we define $g_{a',b'}=e$ or there
are composable edges $a,b$ of $G(X)$ and $\psi _{a'}=\psi _{a},
\psi _{b'}=\psi _{b}$. In this case we define $g_{a',b'}=g_{a,b}$.

One can see using presentations that the fundamental group of
$G(X')$ is isomorphic to the fundamental group of $G(X)$. We
explain this in detail now. It might be helpful for the reader to
draw the barycentric subdivision of a 2-simplex while following
our argument.

Let $T$ be the maximal tree that we pick for the presentation of
the fundamental group of $G(X)$. We will choose a maximal tree for
$G(X')$ that contains $T$. We focus now on the generators and
relators added by a subdivision of a 2-simplex of $X$, $\sigma $.
$\sigma $ in $X$ has three edges which correspond to generators.
After subdivision we obtain 4 new vertex groups and 11 edges (6 of
which are subdivisions of old edges). 4 of the edges correspond
,by definition, to the identity homomorphism from old vertex
groups to the four new vertex groups. To obtain $T'$ We add these
4 edges to $T$ (some of course might already be contained in $T$).
 The relations $\psi _a(g)=a^+ga^-$ for these 4
edges together with $a^+=1,\,for\, a^+\in T'$ show that the new
vertex groups do not add any new generators. Let's call the 4
edges we added to $T$, $a_1,a_2,a_3,a_4$. We remark that 2 more
edges, say $b_1,b_2$ correspond by definition to the identity map
and the relations $ a^+b^+=g_{a,b}(ab)^+$ show that these two new
generators are also trivial (the corresponding $g_{a,b}$'s here
are trivial  as the maps that we compose are identity maps). We
define now a homomorphism from the fundamental group of $G(X)$ to
the fundamental group of $G(X')$ (the presentations given with
respect to $T,T'$ respectively). We focus again on the generators
of the 2-simplex $\sigma $. Vertex groups $G_{\tau }$ are mapped
by the identity map to themselves. Each edge of $\sigma $ is
subdivided in two edges, one of which we added to $T'$. We map
each edge to the edge of the subdivision that we did not add to
$T'$. By the definition of $G(X')$ all relators are satisfied so
we have a homomorphism. It remains to see that it is onto. As we
remarked before 6 of the new edges are trivial in the group. The
other ones can be obtained by successive compositions of the edges
contained in the image (together with edges that are trivial). The
relations  $a^+b^+=g_{a,b}(ab)^+$ for composable edges show that
all generators corresponding to edges are contained in the image.
It is clear that the homomorphism we defined is also 1-1. So it is
an isomorphism.

We return now to our treatment of pairs of splittings.

Let $A_1\star _{C_1}B_1$ (or $A_1\star _{C_1}$), $A_2\star
_{C_2}B_2$ (or $A_2\star _{C_2}$) be a pair of
hyperbolic-hyperbolic splittings of a group $G$ with corresponding
Bass-Serre trees $T_1$, $T_2$. We consider the diagonal action of
$G$ on $Y=T_1\times T_2$ given by
 $$g(t_1,t_2)=(gt_1,gt_2)$$
where $t_1,t_2$ are vertices of, respectively, $T_1$ and $T_2$ and
$g\in G$. We consider the quotient complex of groups in the sense
of Haefliger . If $X$ is the quotient complex $Y/G$ we denote the
quotient complex of groups by $G(X)$.


We give now a detailed description of $G(X)$. We assume for
notational simplicity that the two splittings are $A_1\star
_{C_1}B_1$ and $A_2\star _{C_2}B_2$ (i.e., they are both
amalgamated products). One has similar descriptions in the other
two cases. In the following, if, for example, the splitting along
$C_1$ is an HNN-extension, namely, $G=A_1*_{C_1}$, then one should
just disregard $B_1$, ${\mathcal B}_1$, etc. When there are
essential differences in the HNN- case we will explain the
changes.

 Let $\A _1=T_2/A_1$,
$\B _1=T_2/B_1$, $\A _2=T_1/A_2$, $\B _2=T_1/B_2$,
 $\C _1=T_2/C_1$, $\C _2=T_1/C_2$ be the quotient graphs of the actions
of
$A_1,B_1,C_1$ on $T_2$ and of $A_2,B_2,C_2$ on $T_1$.
Let $A_1(\A _1),B_1(\B _1),C_1(\C _1)$,
$A_2(\A _2),B_2(\B _2),C_2(\C _2) $
be the corresponding Bass-Serre graphs of groups.
We note now that if $e$ is an edge of $X$ which lifts to an edge of
$T_1$
in $Y$ then the
subgraph of the barycentric subdivision of $X$ perpendicular to
the midpoint of this edge
is isomorphic to $\C _1$, and
if we consider it as a graph of groups using the groups assigned to
the vertices and edges
by $G(X)$ then we get a graph of groups isomorphic to  $C_1(\C _1)$.
We identify therefore this
one-dimensional subcomplex of $G(X)$ with $C_1(\C _1)$ and
in a similar way we define
a subcomplex of $G(X)$ isomorphic to $C_2(\C _2) $ and we call
it $C_2(\C _2)$.

We have the following:
\begin{lem}[Van-Kampen theorem]
\label{lem,vankampen}
Let $\C $ be a connected
1-subcomplex of the barycentric subdivision
of $X$ separating locally $X$ in two pieces. We consider $\C$
as a graph of groups
where the groups of the 0 and 1 cells of this graph are induced by
$G(X)$.
Let $C$ be the image of the fundamental group of this graph into
the fundamental group of $G(X)$.
Then the fundamental group of $G(X)$ splits over C.
\end{lem}

\proof

It follows easily from the presentation of the fundamental group
of $G(X)$ given in \cite {H} or in \cite {BH}.
 A detailed
explanation is given in \cite {BH}, ch. III, 3.11 (5), p. 552.
\qed

Since $Y=T_1 \times T_2$ is a product we sometimes use terms
"perpendicular" and "parallel" for certain one-dimensional subsets
in $Y$. Formally speaking, let $p_1,p_2$ be the natural
projections of $Y$ to $T_1,T_2$. Let $e$ be an edge of $T_1$ and
$v$ a vertex of $T_2$. For a point $x \in (e \times v) \subset Y$,
we say that the set $p_1^{-1}(p_1(x))$ is perpendicular to $(e
\times v)$ at $x$.

For convenience we say that the $T_2$-direction is 'vertical' and
the $T_1$ direction is 'horizontal'.

 More formally a set of the form $p_2^{-1}(p_2(x))$ $(x\in Y)$ is
'horizontal'.  We also say that two vertical sets are "parallel".
In the same way, all "horizontal sets", which are of the form
$p_2^{-1}(p_2(x))$, are parallel to each other. Also, those terms
make sense for the quotient $Y/G$ since the action of $G$ is
diagonal, so we may use those terms for the quotient as well.

\begin{defn}[Core subgraph]\label{defn,coregraph}
 Let $A$ be a finitely generated group acting on a tree $T$.
 Let $\A =T/A$ be the quotient graph and let $A(\A )$ be
the corresponding graph of groups. Let $T'$ be a minimal invariant
subtree for the action of $A$ on $T$.
 We define the core of $A(\A )$ to be the subgraph of groups
of $A(\A )$ corresponding to $T'/A$.
\end{defn}

Note that the core of $A(\A )$ is a finite graph. If  $A$ does not
fix a point of $T$ this subgraph is unique. Otherwise it is equal
to a single point whose stabilizer group is $A$. In what follows
we will assume that $C_1,C_2$ are slender groups.
 Therefore the core of $C_1(\C _1)$ (resp. $C_2(\C _2)$) is a circle
unless $C_1$ (resp. $C_2$) acts on $T_2$ (resp. $T_1$) by a dihedral
action in which case
the core is a segment (which might contain more than one edge).

We give now an informal description of the quotient complex of
groups $G(X)$. This description is not used in the sequel but we
hope it will help the reader gain some intuition for $G(X)$.

We consider the graphs of groups $A_1(\A _1),B_1(\B _1),C_1(\C
_1)$. (Disregard $B_1(\B_1)$ if the splitting along $C_1$ is an
HNN-extension. In what follows, this kind of trivial modification
should be made). There are graph morphisms from $\C _1$ to $\A
_1,\B _1$ coming from the inclusion of $C_1$ into $A_1,B_1$. (If
HNN, both morphisms are to $\A_1$). We consider the complex
$[0,1]\times \C _1$. We glue $0\times \C _1$ to $\A _1$ using the
morphism from $\C _1$ to $\A _1$ and $1\times \C _1$ to $\B _1$
using the morphism from $\C _1$ to $\B _1$. The complex we get
this way is equal  to $X$. The vertex groups are the vertex groups
of $A_1(\A _1),B_1(\B _1)$. There are two kinds of edges: the
(vertical) edges of $A_1(\A _1),B_1(\B _1)$ and the (horizontal)
edges of the form $[0,1]\times v$ where $v$ is a vertex of $\C
_1$. The groups of the edges of the first type are given by
$A_1(\A _1),B_1(\B _1)$. The group of an edge $[0,1]\times v$ is
the group of $v$ in $C_1(\C_1)$.

Finally  the group of a 2-cell $[0,1]\times e$ is just the group of $e$
where $e$ is an edge of $\C _1$.

\begin{prop}[Core subcomplex]\label{prop,corecomplex}
 There is a finite subcomplex $Z$ of $X$ such that the fundamental group
of $G(Z)$ is equal to the fundamental group of
$G(X)$.
\end{prop}

\proof We will show that there is a subcomplex $\tilde Z$, of
$T_1\times T_2$ which is invariant under the action of $G$ and
such that the quotient complex of groups corresponding to the
action of $G$ on $\tilde Z$ is finite. We denote this quotient
complex by $G(Z)$. Clearly the fundamental group of $G(Z)$ is
equal to the fundamental group of $G(X)$.

We describe now how can one find such a complex $\tilde Z$. Let
$e=[a,b]$ be an edge of $T_1$ stabilized by $C_1$ and let $p_1:
T_1\times T_2 \to T_1$ the natural projection. $p_1^{-1}(e)$ is
equal to $T_2\times [a,b]$ and $C_1$ acts on $T_2$ leaving
invariant a line $l$, because $C_1$ is slender. $A_1$ acts on
$p_1^{-1}(a)$ and $B_1$ acts on $p_1^{-1}(b)$. Let $S_1,S_2$ be,
respectively, minimal invariant subtrees of $p_1^{-1}(a)$,
$p_1^{-1}(b)$ for these actions. We note that $l\subset S_1,S_2$
since $C_1$ is contained in both $A_1,B_1$. We can take then
$\tilde Z=G(S_1 \cup \{l \times [0,1]\} \cup S_2)$.

The construction is similar if the splitting over $C_1$ is an
HNN-extension ($G=A_1\star _{C_1}$); we simply take $\tilde
Z=G(S_1\cup \{l \times [0,1]\})$ is this case.

Note that $\{ l \times (0,1)\}/C_1$ embeds in $Z$. If the action
of $C_1$ on the line $l$ is not dihedral, then it is an (open)
annulus and if the action is dihedral then it is a rectangle. In
$Z$, some identification may happen at $\{(l \times \{0\}) \cup (l
\times \{1\})\}/C_1$, so that, for example, $\{ l \times
[0,1]\}/C_1$ can be a closed surface in $Z$.

It is easy to see that $Z=\tilde Z/G$ is a finite complex.
Vertices of $Z$ are in 1-1 correspondence with the union of
vertices of $S_1/A_1\cup S_2/B_1$ and the latter set is finite.
Edges of $Z$ correspond to edges of $S_1/A_1\cup S_2/B_1$ and
vertices of $l/C_1$ while 2-cells are in 1-1 correspondence with
edges of $l/C_1$.
\qed

One can define $Z$ also using our previous description of $G(X)$:
\\
We take $Z$ to be the union of the cores of $A_1(\A _1),B_1(\B
_1)$ and $\{  \text{core of }  C_1(\C _1)\times [0,1] \}$. We have
then that the fundamental group of $G(Z)$ is $A_1\star _{C_1}
B_1$. To see this, consider the (vertical) graph perpendicular to
the midpoint of an edge of the form $v\times [0,1]$, ($v\in \C
_1$). The fundamental group of the graph is $C_1$. This graph
separates $Z$ in two pieces. The fundamental groups  of these
pieces are $A_1,B_1$. So from Lemma \ref{lem,vankampen}
(Van-Kampen theorem) we conclude that the fundamental group of
$G(Z)$ is $A_1\star _{C_1} B_1$.

\subsection{Enclosing groups}
\begin{defn}[a set of hyperbolic-hyperbolic splittings]
\label{defn,setofhyp-hyp}
A set,  $I$,
 of splittings of $G$ over slender subgroups is called a set of
hyperbolic-hyperbolic splittings
if for any two splittings in $I$ there is a sequence of splittings in
$I$
of the form $A_i\star _{C_i}B_i$ or  $A_i\star _{C_i}$, $i=1,...,n$,
such that
the first and the last splitting of the sequence are the given
splittings and
any two successive splittings of the sequence are hyperbolic-hyperbolic.
\end{defn}
We remark that a pair of splittings is hyperbolic-hyperbolic (def.
\ref{defn,type}) if and only if the set containing the 2
splittings is a set of hyperbolic-hyperbolic splittings. This
follows from prop. \ref{prop,dual-minimal}.

\begin{defn}[Enclosing graph decomposition]
\label{defn,enclosing}
Let $I$ be a set of hyperbolic-hyperbolic minimal splittings of a group
$G$ along slender groups.
An enclosing group of $I$, denoted by $S(I)$, is a subgroup
in $G$, which is a vertex group of some graph decomposition
of $G$ with the following properties:
\begin{enumerate}
\item
There is a graph of groups decomposition, $\Gamma $, of $G$
with a vertex, $v$, such that
$S(I)$
is the vertex group of $v$, all edges are adjacent to $v$ and
their stabilizers are slender and
peripheral subgroups of $S(I)$ (see below for
the definition of peripheral subgroups).
$S(I)$ contains
conjugates of all the edge groups of
the splittings of $G$
contained in $I$. Each edge of $\Gamma $ gives
a splitting which is elliptic
with respect to all splittings in $I$.
$\Gamma$ is called an enclosing graph decomposition.
 \item  (rigidity)
Suppose $\C'$ is a graph decomposition of $G$ such that any edge
group is slender and gives a splitting which is elliptic-elliptic
to any of the splittings in $I$. Then $S(I)$ is a subgroup of a
conjugate of a vertex group of $\C'$. \item $S(I)$ is an extension
of the (orbifold) fundamental group of a compact $2$-orbifold,
$\Sigma$, by a group $F$ which is a normal subgroup of some edge
group of a splitting contained in $I$. We say that $\Sigma$ is a
base orbifold of $S(I)$, and $F$ is the fiber group. A subgroup of
a group in $S(I)$ which is the induced extension of the (orbifold)
fundamental group of $\partial \Sigma$ by $F$ is called a
peripheral (or boundary) subgroup. We also consider subgroups of
$F$ and of induced extensions of the (finite cyclic) groups of
 the singular points of $\Sigma$, to be peripheral subgroups as well.
\item
Each edge of $\C$ gives a minimal splitting of $G$.
\end{enumerate}

\end{defn}

\begin{rem}
\begin{enumerate}
\item Peripheral subgroups are always proper subgroups of infinite
index of an enclosing group.
\item An enclosing groups
is not slender except when its base $2$-orbifold has a fundamental
group isomorphic to ${\Bbb Z} \times {\Bbb Z}$ or $({\Bbb
Z}_2*{\Bbb Z}_2)\times {\Bbb Z}$,( i.e. the orbifold is  a torus
or an annulus whose two boundary circles are of cone points of
index $2$).
Those two cases are the only tricky ones  that an enclosing group
may have more than one "seifert structure", i.e., the structure of
the extension migth be not unique. For example, ${\Bbb Z}^3$ has
more than one structures of an extension of ${\Bbb Z}^2$ by ${\Bbb
Z}$.
\end{enumerate}
\end{rem}

%

\subsection{Producing enclosing group}
We will show that an enclosing graph decomposition with an
enclosing vertex group $S(I)$ exists (in fact we construct it) for
$I$ given. We start with the simplest case where there are only
two splittings in $I$. As a first step, in the following
proposition using products of trees, we produce a graph
decomposition $\C$ with a vertex group which has properties
(1),(2),(3) in Def \ref{defn,enclosing}.
 Later we will show that
one can also ensure that $\C$ satisfies (4) as well.

One may wonder what happens if we use $A_2,B_2,C_2$ instead of
$A_1,B_1,C_1$ to construct $Z$ in Prop \ref{prop,corecomplex}. In
fact, if both splittings are minimal, then we get the same finite
complex. This is the idea behind the next proposition.

\begin{prop}[Enclosing groups for a pair of splittings]
\label{prop,enclosinggroup}
Let $A_1\star _{C_1}B_1$ (or $A_1\star _{C_1}$) and
$A_2\star _{C_2}B_2$ (or $A_2\star _{C_2}$)
 be a pair of hyperbolic-hyperbolic splittings
of a finitely generated group $G$ over two slender groups $C_1,C_2$.
Suppose that both splittings are minimal.
Then there is a graph decomposition of $G$ with
a vertex group which has the
properties 1, 2, 3 in Def \ref{defn,enclosing} for $C_1,C_2$.
\end{prop}

\begin{rem}\label{rem,minimal}
If $G$ does not split over a subgroup of infinite index in
$C_1,C_2$ then the two splittings along $C_1,C_2$ are minimal (see
remark \ref{rem,infiniteindex}). This hypothesis is used in \cite
{RS} and \cite {DS} instead of minimality.
\end{rem}

We first construct a graph decomposition of $G$ and show that it
is the desired one for Prop \ref{prop,enclosinggroup} later. Let
$T_i$ be the Bass-Serre tree of the splitting over $C_i$. Consider
the diagonal action of $G$ on $Y=T_1\times T_2$. Let $G(X)$ be the
quotient complex in the sense of Haefliger. Consider the finite
subcomplex of $G(X)$, $G(Z)$, constructed in proposition
\ref{prop,corecomplex}. Let $e$ be an edge of $A_1(\A _1)$ lying
in the core of $C_1(\Gamma _1)\times [0,1]$. In other words $e \in
\A_1 \cap \{\Gamma_1 \times \{0\}\}$. Consider the (horizontal)
graph in $G(Z)$ perpendicular to $e$ at its midpoint. In other
words this is the maximal connected graph passing through the
midpoint of $e$ whose lift to $T_1\times T_2$ is parallel to
$T_1$.
The fundamental group of the graph of groups
corresponding to this graph  is a subgroup of a conjugate of $C_2$.
Since both the assumptions and the conclusions of Prop \ref{prop,enclosinggroup}
do not change if we take conjugates in $G$ of the splittings
along $C_1,C_2$,  without loss
of generality (by substituting the splitting along $C_2$ by a conjugate
in $G$) we can assume
that this graph is a subgraph of $C_2(\C _2)$.

\noindent {\it claim}. Consider the squares intersecting the core
of $C_2(\C_2) $. Then, this set of squares contains all squares of
$G(Z)$.

\noindent To argue by contradiction, we distinguish two cases
regarding the set $Z \cap (T_1/C_2 \times \{1/2\})$. We naturally
identify the core of $C_2(\Gamma_2)$ with a subgraph in $T_1/C_2
\times \{1/2\}$.

\noindent
(i) If $Z$ does not contain the core
of $C_2(\C _2)$ then, by Lemma \ref{lem,vankampen}(Van-Kampen),
$G$ splits over an infinite
index subgroup of
 $C_2$. Moreover,
  $A_1\star _{C_1}B_1$ (or $A_1*_{C_1}$) is
hyperbolic-elliptic with respect to this new splitting
contradicting our hypothesis. We explain this in more detail. The
type of the splitting over $C_1$, i.e., either an amalgamation or
HNN extension, does not make difference in this discussion. On the
other hand, we may need to make a minor change according to the
type of the splitting along $C_2$, which we pay attention to. What
requires  more attention, because the topology of the base
$2$-orbifold of $S(I)$ becomes different, is  the type of the
action of $C_2$ on $T_1$, i.e., dihedral or not, although the type
of the action of $C_1$ on $T_2$ is not important once the
subcomplex $Z$ is constructed.

Let's first assume that  the splitting along $C_2$ is not dihedral.
 Let $l_2 \subset T_1$ be the invariant
line of the action by $C_2$. The core of $C_2(\C_2)$ is $l_2/C_2$,
which is a circle by the assumption we made. The circle $l_2/C_2$
is a retract of a graph $T_1/C_2$, so that $T_1/C_2 - l_2/C_2$ is
a forest, i.e., each connected component is a tree. Therefore, if
$Z$ does not contain the core of $T_1/C_2 \times \{1/2\}$, then $Z
\cap (T_1/C_2 \times \{1/2\})$ is a forest. Let $U_1, \cdots, U_n$
be the connected components of the forest. Then, if we cut $Z$
along each $U_i$, and apply Lemma \ref{lem,vankampen}, we get a
splitting of $G$ along the fundamental group (in the sense of
graph of groups) of $U_i$, $K_i$, which is a subgroup of infinite
index in $C_2$.

By construction, $K_i$ is contained in $C_1$. Therefore this
splitting along $K_i$ is elliptic with respect to $A_1*_{C_1}B_1$
(or $A_1*_{C_1})$. Moreover, $C_1$ is hyperbolic to at least one
of the splittings along $K_i$'s. This is because if not, then
$C_1$ is contained in a conjugate of $A_2$ or $B_2$ (or $A_2$ in
the case that the splitting along $C_2$ is an HNN-extension),
which is impossible, since $C_1$ is hyperbolic with respect to
$A_2*_{C_2}B_2 ({\rm or }A_2*_{C_2})$. The last claim does not
require the theory of complex of groups, but just Bass-Serre
theory; since each $U_i$ is a tree, $U_i$ contains a vertex,
$u_i$, whose vertex group is $K_i$. Let $e_i=u_i \times [0,1]
\subset Z$. If we delete all (open) squares and edges parallel to
$e_i$ (except for $e_i$) in $U_i \times [0,1]$ from $Z$, the
fundamental group (in the sense of Haefliger) does not change, and
also the edge $e_i$ gives the splitting of $G$ along $K_i$. If we
do the same thing for all $U_i$'s, the subcomplex of $Z$ we obtain
is indeed a graph, $\Lambda$, where the edges $e_i$'s are parallel
to each other, and no other edges are parallel to them. Note that
all of those other edges are the ones which were in the graph
$\A_2 \cup \B_2 \subset Z$ (or just $\A_2$ in the case of
HNN-extension). Therefore, if $C_1$ is elliptic with respect to
all the splittings along $K_i$, it means that $C_1$ is conjugated
to the fundamental group (in the sense of Bass-Serre) of a
connected component of $\Lambda -\cup_i U_i$, which is a subgraph
of either $\A_2$ or $\B_2$ (or just $\A_2$ in the case of
HNN-extension). This means that $C_1$ is a conjugate of a subgroup
of either $A_2$ or $B_2$ (or $A_2$ in the case of HNN-extension).
This is what we want.

We are still left with the case that $C_2$ is dihedral on $T_1$.
The argument only requires a notational change. $T_1/C_2 \times
\{1/2\}$ is a forest and we look at each connected component, and
appropriately delete all squares and some edges from $Z$ without
changing the fundamental group, which is $G$, and get a graph of
groups at the end as before. We omit details.

We remark that our argument does not change if the action of $C_1$ on $T_2$ is dihedral or
not. So we treated all possibilities in terms of the type of the splittings along $C_1,C_2$ and also
the type of the actions of $C_1,C_2$.

\noindent (ii) If on the other hand $Z \cap (T_1/C_2 \times
\{1/2\})$ is bigger than the core of $C_2(\C _2)$ we can delete
from $G(Z)$ the $2$-cells (i.e., squares) containing edges of this
graph which do not belong to the core of $C_2(\C _2)$ without
altering the fundamental group. To explain the reason, let's first
suppose that the action of $C_2$ on $T_1$ is not dihedral. Then
the core is topologically a circle, $c_2$. The connected
component, $U$, of the finite graph $Z \cap (T_1/C_2 \times
\{1/2\})$ which contains the core is topologically the circle with
some trees attached. The fundamental group (in the sense of graph
of groups) of not only the circle $c_2$ , but also the graph $U$
is $C_2$.Therefore, one can remove those trees from $U$ without
changing the fundamental group, which is $C_2$.

In $Z$, $U \times (0,1)$ embeds, and one can remove the part
$(U-c_2) \times (0,1)$ from $Z$ without changing its fundamental
group. One can see this using the presentation of the fundamental
group of $G(Z)$. For the reader's convenience we give also an
argument using the action of $G$ on $\tilde Z$. Let $p_2:\tilde
Z\to T_2$ be the natural projection from $\tilde Z$ to $T_2$. if
$e$ is an edge of $T_2$, $p_2^{-1}(e)$ is a connected set of the
form $L_e\times e$. Let
 $Stab(e)$ be the stabilizer of $e$ in $T_2$ (which is a
conjugate of $C_2$) and $l_e$ the line invariant under $Stab(e)$
on $T_1$. Then by the discussion above $L_e$ contains $l_e$ and is
connected. We will show that $L_e=l_e$. Indeed if not we consider
the subcomplex of $\tilde Z$ obtained by the union of $l_e\times
e$ over all edges $e\in T_2$ with $p_2^{-1}(v)$ over all vertices
$v\in T_2$. Let's call this complex $\tilde Z_1$. It is clear that
$\tilde Z_1$ is connected, simply connected and invariant under
the action of $G$. The quotient complex of groups $G(Z_1)$ is a
subcomplex of $G$. If for some $e$, $L_e\ne l_e$ $G(Z_1)$ is
properly contained in $G(Z)$. By the preceding discussion it
follows that the splittings over $C_2,C_1$ are
hyperbolic-elliptic, a contradiction.

In the case when the action of $C_2$ on $T_1$ is dihedral, then
the core is a segment, and $Z \cap (T_1/C_2 \times \{1/2\})$ is a
graph which is the segment with some finite trees attached. In
this case one can delete the squares which contain those trees
from $Z$ without changing the fundamental group of $Z$, which is
$G$.

But then the complex obtained, after deleting those unnecessary squares,  does not
contain the core of $C_1(\C _1)$, because there are no other squares in $Z$
than the ones which contains the core of $C_1(\C_1)$, which implies
that $G$ splits over an infinite index subgroup of $C_1$,
and this new splitting
is elliptic-hyperbolic with respect to $A_2\star _{C_2}B_2$
 (or $A_2*_{C_2}$), which is a
contradiction. The last part follows from the same consideration
as the last part of the case (i), so we omit the details. We
showed the claim.

From this claim it follows that if we apply the same construction
of $Z$ in Prop \ref{prop,corecomplex} using $A_2,B_2$ instead of
$A_1,B_1$, the resulting complex contains the same set of squares.

We describe the topology of $Z$. If none of $C_1,C_2$ acts as a
dihedral group
 the above implies that every edge in $Z$ which is a side of a $2$-cell
lies on exactly two $2$-cells. Therefore the link of every vertex
of $Z$ is a union of disjoint circles and points. It then follows
that the union of $2$-cells in $Z$ is topologically a closed
surface with, possibly, some (vertex) points identified. $Z$ is
this $2$-dimensional object with some graphs attached at vertices;
if one deletes from $Z$ those graphs including attaching vertices
and identified vertices, one obtains a compact surface with
punctures.

 If at  least
one of $C_1,C_2$ acts as a dihedral group then $Z$ is
topologically a compact surface with boundary with, possibly, some
points identified and some graph attached. Therefore the links of
vertex points on this surface are disjoint unions of circles,
segments and points. The boundary components come from the
dihedral action(s), and there are at most 4 connected components.
To see this, suppose that only $C_1$ is dihedral on $T_2$, and let
$l_1$ be its invariant line. Then the rectangle $l_1/C_1 \times
(0,1)$ embeds in the surface. Let $u_1,u_2$ be the boundary points
of the segment $l_1/C_1$. Then the edges $u_1 \times [0,1], u_2
\times [0,1]$ are exactly the boundary of the surface. Note that
$u_1 \times (0,1)$ embeds, but possibly, $u_1 \times [0,1]$ may
become a circle in $Z$. $u_1 \times [0,1]$ and $u_2 \times [0,1]$
may become one circle in $Z$ as well. Therefore, the surface has
at most two boundary components in this case. If $C_2$ is dihedral
on $T_1$ as well, then there are two more edges which are on the
boundary, so that there are at most 4 boundary components.

\begin{rem}[Priority among splittings]\label{rem,priority}
Let $Z$ be the complex we constructed  in the proof of Prop.
\ref{prop,corecomplex}. Let $l_1$ be the invariant line in $T_2$
by $C_1$ and $c_1=l_1/C_1$. We may call $c_1$ the core of $C_1$.
If the action of $C_1$ is dihedral, then $c_1$ is a segment, or
else, a circle. The core $c_1$ embeds in $Z$, and if we cut $Z$
along $c_1$ we get (not a conjugate, but exactly) the splitting
$A_1*_{C_1} B_1$ (or $A_1*_{C_1}$). Similarly, let $l_2$ be the
invariant line in $T_1$ by $C_2$, and $c_2=l_2/C_2$. As before,
this core $c_2$ is either a segment or a circle, and embeds in
$Z$. Cutting $Z$ along $c_2$ , we get a splitting of $G$ along a
conjugate of $C_2$. But, unlike the splitting along $C_1$, this
splitting may be different from the original splitting along
$C_2$. This point becomes important later, that we can keep at
least one splitting unchanged (along $C_1$ in this case) in $Z$,
because we gave priority to the splitting along $C_1$ over $C_2$
when we constructed $Z$.

However, it is true that if $G$ does not split along a subgroup in
$C_1$ of infinite index, then the new splitting along $C_2$
obtained by cutting $Z$ along $c_2$ is the same (up to
conjugation) as the original one. It is because that under this
assumption, $Z$ does not have any graphs attached, and it is just
a squared complex.

Although the new splitting along $C_2$ may be different from the original one,
it is hyperbolic-hyperbolic with respect to the splitting along $C_1$.
It follows from Lemma \ref{prop,dual-minimal} that the new splitting
along $C_2$ is minimal.

\end{rem}
We now explain how to obtain the desired graph of groups
decomposition of $G$. First, let the group $S=S(C_1,C_2)$ be the
subgroup of the fundamental group of
 $G(Z)$ corresponding to the subcomplex of $G(Z)$ which is the union of
the cores of $C_1(\C _1)$ and $C_2(\C _2)$, namely, $S$ is the
image in $G$ of the fundamental group (in the sense of a graph of
groups) of this union. Here we use Haefliger's notation; the cores
of $\C _1,\C _2$ are contained in the barycentric subdivision of
$Z$ which is used in the definition of its fundamental group.

Using lemma \ref{lem,vankampen}(Van-Kampen theorem) we show that
$G(Z)$ is the fundamental group of a graph of groups, which we
call $\C$. The vertices of this graph are as follows: there is a
vertex for each connected  component of $Z$ minus the cores of
$C_1(\C _1)$ and $C_2(\C _2)$.
 The vertex group is
the fundamental group of the component in Haefliger's sense. We
remark that each such component contains exactly one vertex group
of the 'surface' piece of $Z$ with, possibly, a graph attached at
the vertex. The fundamental group of the component is then the
fundamental group of the graph of groups of the attached graph
(and is equal to the group of the vertex if there is no graph
attached).

 There is also a vertex with group $S(C_1,C_2)$. There is an edge
for each component of the intersection
 between the union of the cores of $C_1(\C _1)$, $C_2(\C  _2)$ and
each vertex component.

Note that such an intersection is topologically a circle or a segment
(this happens only when at least one of the actions of $C_i$ is
dihedral).
  As this intersection is a subgraph of
the union of the cores of $C_1(\C _1)$, $C_2(\C _2)$ there is a
group associated to it, namely the image of the fundamental group
of this subgraph in $G(Z)$. Note that the graph of groups that we
described here is a graph of groups in a generalized sense, i.e.,
the edge groups do not necessarily inject into the vertex groups.
Note that every vertex group except $S$ injects in $G$.

To understand the group $S$ better, let $U$ be the union of the
cores of $C_1(\C _1)$, $C_2(\C _2)$, which is a graph in $Z$. If
we consider a small closed neighborhood, $\bar U$, of $U$ in $Z$,
it is a compact surface with boundary in general. We may consider
the graph, $P$, corresponding to an edge, $e$, of $\Gamma$ as a
subset in the boundary of $\bar U$, which is either a circle or a
segment. Let $F<G$ be the stabilizer of a square in $Z$. Then, the
fundamental group, in the sense of Haefliger, of $P$ is an
extension of ${\Bbb Z}$ (when $P$ is a circle) or ${\Bbb Z}_2
* {\Bbb Z}_2$ (when $P$ is a segment) by $F$.
Since the group $F$ is a subgroup of $G$, the image of the
fundamental group of $P$ in $G$ is an extension
of (1) ${\Bbb Z}$, (1') ${\Bbb Z}_n$, (1'') the trivial group;
(2) ${\Bbb Z}_2*{\Bbb Z}_2$, (2') ${\Bbb Z}_2$, or
(2'') a finite dihedral group of order $2n$ by $F$.
As a consequence, the image of the fundamental group of $U$
(as well as $\bar U$) in $G$, which is $S$ by our definition, is an extension of the
orbifold fundamental group of a $2$-orbifold, $\Sigma$, by $F$ such that
$\Sigma$ is obtained from the compact surface $\bar U$ by
adding to each $P$  (1) nothing, (1') a disk with a cone point of index $n$
at the center, (1'') a disk; (2) nothing,
(2') a half disk such that the diameter consists of cone points
of index $2$ (in other words, we just collapse the segment $P$ to a point),
or (2'') a half disk such that the diameter consists of
cone points of index $2$ except for the center whose index is $2n$.
The orbifold fundamental group of this
$2$-orbifold is $S$. Note that $\Sigma$ is no more embedded in $Z$,
but the surface $\bar U$ is a subsurface of $\Sigma$.

\begin{rem}
By our construction of $Z, U$,  $\Sigma$, there is a simple closed
curve or a segment (the core) on $\Sigma$ which corresponds to
each of $C_1,C_2$. If we cut $\Sigma$ along it, we obtain a
splitting of $S$ along $C_1$ or $C_2$, respectively, which also
gives a splitting of $G$, as we do by cutting $Z$. Although one of
them may be different from the original one, both of them are
minimal (use Prop \ref{prop,dual-minimal}).
\end{rem}

\subsection{Proof of Prop \ref{prop,enclosinggroup}}

\proof We will show that the graph decomposition $\C$ with
$S(C_1,C_2)$ we constructed satisfies the properties 1,2,3  in Def
\ref{defn,enclosing}.
In fact $S(C_1,C_2)$ is an enclosing group
for $C_1,C_2$ although we may need to modify $\C$ so that the
property (4) holds as well. We will discuss this point later.

(3) is clear. By construction, $S$ is an extension of the
orbifold fundamental group of the compact $2$-orbifold $\Sigma$
by a group $F$, which is the edge stabilizer subgroup in $C_1$
when it acts on the tree $T_2$, hence a normal subgroup of $C_1$.
$F$ is slender since it is a subgroup of a slender group $C_1$.

(1). Let $v$ be the vertex of $\Gamma$ whose
vertex group is $S$. By our construction, all edges
are adjacent to $v$. The edge group, $E$, of an edge is slender
since there is the following exact sequence;
$1 \to F \to E \to Z \to 1$ such that the group $Z$ is
either the trivial group, the fundamental group of one
of the singular points of $\Sigma$, so that isomorphic to ${\Bbb Z}_n$,
or a subgroup of the (orbifold) fundamental group
of $\partial \Sigma$, so that isomorphic to ${\Bbb Z}_2
*{\Bbb Z}_2$. In any case, $E$ is slender and a peripheral subgroup
of $S$.
Clearly $S$ contains conjugates of $C_1,C_2$, because
$\Sigma$ contains the graph $U$, which is the union
of the cores for $C_1,C_2$.
By construction of $\Gamma$, all vertex groups except for $S$ are elliptic
on both $T_1$ and $T_2$, so that all edge groups of $\Gamma$
are elliptic on $T_1,T_2$ since they are subgroups of vertex groups.
We showed (1).


To prove that enclosing groups are `rigid', namely
the property (2) in Def \ref{defn,enclosing}, we recall some
results from Bass-Serre theory.
\begin{prop}[Cor 2 in \S 6.5 of \cite{Se}]\label{prop,bassserre}
Suppose $G$ acts on a tree. Assume $G$ is generated
by $s_1,\cdots, s_l$ and all $s_i$ and $s_i s_j (i \not= j)$
are elliptic on the tree. Then $G$ is elliptic.
\end{prop}

Let $c$ be a simple closed curve on $Z$ which avoids vertices of
$Z$. Using repeated barycentric subdivisions of $G(Z)$ we see that
$c$ is homotopic in $Z-Z^{(0)}$ to a curve lying in the 1-skeleton
of the iterated barycentric subdivision. Let's assume then that
$c$ is a curve lying in the 1-skeleton of an iterated barycentric
subdivision. Then cutting $Z$ along $c$, we get a splitting of $G$
along the fundamental group (in the sense of graph of groups) of
$c$, (\ref{lem,vankampen}). If the fundamental group of $c$ is not
contained to a vertex group, then the splitting induced by $c$ is
hyperbolic-hyperbolic with respect to either the splitting along
$C_1$ or $C_2$, so that in particular, it is non trivial and a
minimal splitting. We call such simple closed curves $c$ {\it
essential }

In the case $Z$ has a boundary (i.e., there exists an edge which
is contained in only one square), let $c$ be an embedded segment
whose boundary points are in the boundary of $Z$. Cutting $Z$
along $c$, one also obtains a splitting of $G$ along the
fundamental group of $c$.  If the fundamental group of $c$ is not
contained to the fundamental group of $\partial Z$ then this
splitting is hyperbolic-hyperbolic with respect to at least one of
the splittings along $C_1$ and $C_2$, so it is non trivial and
minimal. We call such segments $c$ {\it essential}. We remark that
essential simple closed curves and essential segments correspond
to subgroups of the fundamental group of $\Sigma $.

If $\partial \Sigma$ contains segments of singular points of index
two (reflection points) we denote this set by $(\partial \Sigma;2)$.

\begin{cor}\label{cor,elliptic}
If all splittings over the slender groups which are represented by
essential simple closed curves on $\Sigma $ and essential embedded
segments are  elliptic on $\C$, then $S=S(C_1,C_2)$ is elliptic on
$\C$.
\end{cor}

\proof We explain how to choose a finite set of generators of $S$
so that we can apply Prop \ref{prop,bassserre}. First choose a
finite set of generators $f_i$ of $F$ ($F$ is slender, so that
finitely generated).

 Let's
first assume that $(\partial \Sigma;2)=\emptyset$. Then, one can
choose a set of non-boundary simple closed curves $c_1, \cdots,
c_l$ on $\Sigma$ so that all $c_i c_j (i \not= j)$ are also
represented by simple closed curves and that the elements
corresponding to $c_i$ generate the fundamental group of $\Sigma$.
Each $c_i$ or $c_i c_j (i \not= j)$ represents a slender subgroup
in $G$ with the fiber group $F$, which gives a splitting of $S$.
By assumption all of those splittings are elliptic on $\C$.
Therefore we apply Prop \ref{prop,bassserre} to $S$ with the
generating set of $\{ f_i,c_j \}$ and conclude that $S$ is
elliptic on  $\C$.

In the case $(\partial \Sigma;2) \not =\emptyset$, we need extra
elements represented by embedded segments $(s, \partial s) \subset
(\Sigma, (\partial \Sigma;2))$. Put an order to the connected
components of $(\Sigma, (\partial \Sigma;2))$, and take a finite
set of embedded segments, $s_i$, such that any adjacent (in the
order) pair of components of $(\partial \Sigma;2)$ is joined by a
segment. Then the set $\{ f_i, c_j, s_k \}$ generates a subgroup
$S'<S$ of finite index. By Prop \ref{prop,bassserre}, $S'$ is
elliptic on $\C$, so that so is $S$. \qed(Cor \ref{cor,elliptic}).

We now show that enclosing groups are `rigid'.
\begin{lem}[Rigidity]\label{lem,rigidity}
Let
$A_1\star _{C_1}B_1$ (or $A_1\star _{C_1}$) and $A_2\star _{C_2}B_2$
(or $A_2\star _{C_2}$)
be as in proposition \ref{prop,enclosinggroup} and
 let $S=S(C_1,C_2)$ be the group
constructed above.
Suppose $\C'$ is a graph decomposition of $G$ such that
any edge group is slender and elliptic
to both of the splittings over $C_1,C_2$.
Then $S$ is a subgroup of a conjugate of a vertex group of $\C'$.
\end{lem}

\proof $\Gamma$ denotes the graph of groups decomposition we
constructed with $S$ as a vertex group. As we pointed out in
Remark \ref{rem,priority}, although the splitting of $G$ over
$C_2$ which we obtain by cutting $\Sigma$ along the core curve for
$C_2$ may be different from the original one, this splitting is
minimal, because it is hyperbolic-hyperbolic to the (original)
splitting over $C_1$, (see Prop. \ref{prop,dual-minimal}).

Let $T,T'$ be the Bass-Serre trees of $\Gamma, \Gamma '$. Our goal
is  to show that $S$ is elliptic on $T'$. Let $c \subset \Sigma$
be an essential simple closed curve  or $(s, \partial s) \subset
(\Sigma, (\partial \Sigma;2))$ an embedded essential segment, and
$C <S$ the group represented by it. If we show that $C$ is
elliptic on $T'$, then Prop 2.9 implies that $S$ is elliptic on
$T'$. The splitting of $G$ over $C$ by cutting $\Sigma$ along $c$
or $s$  is minimal by Prop. \ref{prop,dual-minimal} since it is
hyperbolic-hyperbolic to one of the minimal splittings over
$C_1,C_2$.

Let $e$ be an edge of $\C'$ with edge group, $E$. Since the
subgroup $E$ is elliptic with respect to the splittings over
$C_1,C_2$ (i.e., elliptic on the both trees for the two
splittings), it fixes a vertex of $T_1\times T_2$. Therefore it is
contained in a conjugate of a vertex group of $G(Z)$, which is not
$S$. It follows that the group $E$ is elliptic  with respect to
the splitting of $G$ over $C$. Since the splitting along $C$ is
minimal, by Prop \ref{prop,dual-minimal}, it is elliptic-elliptic
with respect to the splitting of $G$ over $E$ obtained from $\C'$
by collapsing all edges but $e$. Since $e$ was arbitrary, the
subgroup $C$ is elliptic on $T'$.
%
\qed(lemma \ref{lem,rigidity}).

Lemma \ref{lem,rigidity}
implies (2) in Def \ref{defn,enclosing}.
We have verified the items (1),(2),(3) in Def \ref{defn,enclosing}
 for $S(C_1,C_2)$
which finishes the proof of Prop \ref{prop,enclosinggroup}.
\qed(Prop \ref{prop,enclosinggroup}).

\begin{rem}[Maximal peripheral subgroup]
\label{rem,maxperipheral} Note that the subset of $\partial
\Sigma$ which is produced by the cutting of $Z$ is exactly
$\partial \Sigma - {\rm interior \, of}(\partial \Sigma;2)$. Let
$c$ be a connected component of this subset, and $E$ the
corresponding (peripheral) subgroup of $S$. Let's call such
peripheral subgroup of $S$ {\it maximal}. $E$ is an edge group of
$\Gamma$ by our construction, so that any maximal peripheral
subgroup of $S$ is an edge group of $\Gamma$. For example when the
fiber group $F$ is trivial, $\Sigma$ is a $2$-manifold with
boundary. Then the infinite cyclic subgroup in $S=\pi_1(\Sigma)$
corresponding to each boundary component of $\Sigma$ is an edge
group of $\Gamma$. In this sense, $\Sigma$ does not have any free
boundary points.
\end{rem}

\subsection{Producing enclosing graph decomposition}
We now discuss the property (4)  of Def \ref{defn,enclosing}.
In general
the edges of $\C$ we obtained in Prop \ref{prop,enclosinggroup}
 may give non-minimal
splittings. See the example.

However, by applying Prop \ref{prop,modification} to $\C$, there
is a refinement (see Def \ref{defn,refinement}) of $\C$ such that
all edges give minimal splittings of $G$. We then verify that the
refinement satisfies all the properties of Def
\ref{defn,enclosing}, most importantly, $S$ remains a vertex
group, and is the enclosing group of the decomposition we get.

\begin{example*}
This example is suggested by V.Guirardel to us. We thank him.
Let $G={\Bbb Z}^3 *A$ such that $A$ is a non-trivial group.
Fix free abelian generators $a_1,a_2,a_3$ of ${\Bbb Z}^3$.
Write ${\Bbb Z}^3$ as an HNN-extension 
${\Bbb Z}^2*_{{\Bbb Z}^2}$ such that ${\Bbb Z}^2=\langle a_2,a_3 \rangle$
and the stable letter is $a_1$. This extends to an HNN-extension 
$G=({\Bbb Z}^2*A)*_{{\Bbb Z}^2}$ over ${\Bbb Z}^2=\langle a_2,a_3 \rangle$.
Let $T_1$ be the Bass-Serre tree of this splitting. We abuse the notation 
and call the splitting $T_1$ as well.
Similarly, we obtain HNN-extensions $T_2$ and $T_3$:
$G=({\Bbb Z}^2*A)*_{{\Bbb Z}^2}$ over ${\Bbb Z}^2=\langle a_1,a_3 \rangle$, 
and $\langle a_1,a_2 \rangle$, with Bass-Serre trees $T_2,T_3$.
For $i \not= j$, the pair of splittings $T_i,T_j$ is 
hyperbolic-hyperbolic.
Each splitting $T_i$ is minimal, because if $G=({\Bbb Z}^2*A)*_{{\Bbb Z}^2}$ 
was not minimal, then 
the corresponding HNN-extension ${\Bbb Z}^3={\Bbb Z}^2*_{{\Bbb Z}^2}$
would give (for example, by taking product of trees)
a splitting of ${\Bbb Z}^3$ over 
${\Bbb Z}$ or the trivial group, which is impossible. 
We obtain a graph decomposition for the pair $T_1,T_2$
by taking product of trees: $G={\Bbb Z}^3*_{{\Bbb Z}}({\Bbb Z}*A)$
such that this is an amalgamation over $\langle a_3 \rangle$
with two vertex groups ${\Bbb Z}^3$ and $\langle a_3 \rangle *A$.
Let's call this decomposition, and its Bass-Serre tree $T$.
The vertex group ${\Bbb Z}^3$ is the enclosing group such that 
the base is a torus with the fundamental group $\langle a_1,a_2 \rangle$
and the fiber group is $\langle a_3 \rangle$.
This splitting along $\langle a_3 \rangle$ is not minimal, 
because it is hyperbolic-elliptic to $T_3$.
We now demonstrate how to handle this problem 
using Prop \ref{prop,modification}.
Following the first proof of Prop \ref{prop,modification},
we refine $T$ such that all edge gives a minimal splitting.
Take product of trees of $T,T_3$.
The core is topologically an annulus, which 
contains one square, and three edges, where two of them are 
vertical, and they are loops. Since $A$ fixes a vertex of $T \times T_3$,
$A$ is a vertex group of the core.
We remove the (open) square, and also one vertical loop,
whose edge group is trivial, 
appropriately without changing the fundamental group.
We obtain a graph decomposition with two edges, 
$G=A* {\Bbb Z}^2* _{{\Bbb Z}^2}$ such that both ${\Bbb Z}^2$ are 
$\langle a_1, a_2 \rangle$.
Next, collapse the other vertical loop, whose edge group 
is $\langle a_1,a_2 \rangle$, in the graph decomposition.
We are left with the horizontal edge, whose edge group is trivial,
and obtain $G={\Bbb Z}^3 *A$, which is a refinement of
$T$. This is an enclosing decomposition for $T_1,T_2$ with the enclosing 
vertex group ${\Bbb Z}^3$.
\end{example*}

\begin{lem}[Refinement of $\C$]\label{lem,edgerefinement}
There is a refinement, $\C'$, of $\C$
such that
\begin{enumerate}
\item
each edge of $\C'$ gives a minimal splitting of $G$,
\item
each edge group of $\Gamma'$ is a subgroup of some
edge group of $\Gamma$,
\item
$S=S(C_1,C_2)$ remains a vertex group of $\C'$,
\item
each edge group is a peripheral subgroup of $S$.
\end{enumerate}
\end{lem}

\proof
Apply Prop \ref{prop,modification}
 to $\C$ and obtain a refinement $\C'$ such
that each edge of $\Gamma'$ gives a minimal splitting of $G$. Each
edge group, $E$, of $\C'$ is a  subgroup of some edge group of
$\C$. Therefore, by Prop \ref{prop,enclosinggroup},
 $E$ is  elliptic to both splittings of $G$
along $C_1,C_2$.
Since both of the splitting
along $C_1,C_2$ are  minimal, each of them
is  elliptic-elliptic with respect to the
splitting over $E$.
Therefore by Lemma \ref{lem,rigidity},
$S=S(C_1,C_2)$ is a subgroup of a conjugate
of
 some vertex group, $V$, of $\C'$. But since $\C'$ is a refinement of
$\C$
and $S$ is a vertex
group of $\C$, $S$ is a conjugate of $V$.
$E$ is a peripheral subgroup of $S$ since it is a
subgroup of a peripheral subgroup.
\qed

We collapse all edges in $\C'$ which are not adjacent to the vertex
whose vertex group is $S$, and still call it $\C'$.
Then by Prop \ref{prop,enclosinggroup} and
Lemma \ref{lem,edgerefinement},
$\Gamma'$ is an enclosing graph decomposition
with an enclosing vertex group $S(C_1,C_2)$
for the splittings along $C_1,C_2$.
We have shown the following.

\begin{prop}[Enclosing decomposition for a pair of splittings]
\label{prop,enclosingdecomposition}
Let $A_1\star _{C_1}B_1$ (or $A_1\star _{C_1}$) and
$A_2\star _{C_2}B_2$ (or $A_2\star _{C_2}$)
 be a pair of hyperbolic-hyperbolic splittings
of a finitely generated group $G$ over two slender groups $C_1,C_2$.
Suppose that both splittings are minimal.
Then an enclosing graph decomposition of $G$
exists for those two splittings.
\end{prop}

\section{JSJ-decomposition}

\subsection{Dealing with a third splitting}
Let $G$ be a finitely presented group.
We want to show that an enclosing graph decomposition exists
for a set, $I$, of hyperbolic-hyperbolic minimal splittings
of $G$ along slender groups.
We already know this when $I$ contains only two elements
by Prop \ref{prop,enclosingdecomposition}. We now discuss
the case when $I$ has three elements.

\begin{prop}[Enclosing group]
\label{prop,enclosing}
Let $I$ be a set of hyperbolic-hyperbolic splittings (Def \ref{defn,setofhyp-hyp})
of a finitely generated group $G$. Suppose
all of them are minimal splittings. Suppose
that $I$ consists of three splittings.
Then an enclosing graph decomposition
 of $G$ exists for $I$.
\end{prop}

\proof Suppose that the three splittings in $I$ are along
$C_1,C_2,C_3$. We may assume that the pair of the splittings along
$C_1,C_2$, and also the pair for $C_2,C_3$ are
hyperbolic-hyperbolic. Apply Prop \ref{prop,enclosinggroup} to the
first pair, and obtain an enclosing graph decomposition, $\C$,
with the vertex group $S=S(C_1,C_2)$. We remark that $S(C_1,C_2)$
depends on the two splittings, not only the two subgroups. Note
that by cutting the $2$-orbifold $\Sigma$ for $S$ along a simple
closed curve or a segment corresponding to each of $C_1,C_2$, we
obtain a minimal splitting of $G$ along $C_1$, and $C_2$,
respectively. Although this splitting along $C_2$ is possibly
different from the original one, it is still a minimal splitting
so it is hyperbolic-hyperbolic with respect to the splitting along
$C_3$.

Let's assume first that the group $C_3$ is elliptic with respect
to $\Gamma$. Then $C_3$ is a subgroup of a conjugate of $S$. This
is because if $C_3$ was a subgroup of a conjugate of a vertex
group of $\Gamma$ which is not $S$, then the group $C_3$ is
elliptic with respect to both of the (new) splittings of $G$ along
$C_1,C_2$ which we obtain by cutting $\Sigma$. It then follows
that the splitting along $C_3$ would be elliptic-elliptic with
respect to both of the original splittings of $G$ along $C_1,C_2$,
which is a contradiction. Let $\Gamma'$ be a refinement of
$\Gamma$ which we obtain by Prop
\ref{prop,enclosingdecomposition}, which is an enclosing graph
decomposition for the splittings along $C_1,C_2$. We claim that
$\Gamma'$ with an enclosing vertex group $S$ is an enclosing
decomposition for the three splittings. First, the properties
2,3,4 are clear. To verify (the non-trivial part of) the property
1, let $e$ be an edge of $\Gamma '$ with edge group, $E$. We want
to show that the group $E$ is elliptic with respect to the
splitting in $I$ along, $C_3$. We know that $C_3<S$ by our
assumption. Since the group $S$ is elliptic with respect to the
splitting of $G$ along $E$ which the edge $e$ gives, so is $C_3$.
Since both of the splittings along $C_3$ and $E$ are minimal, it
follows that the group $E$ is elliptic with respect to the
splitting along $C_3$. This proves property 1.

We treat now the case that $C_3$ is hyperbolic with respect to
$\Gamma$. This is the essential case. Let $T_{\C},T_3$ be,
respectively, the Bass-Serre trees of $\C $ and the splitting over
$C_3$. Since the splitting along $C_3$ is minimal, there is at
least one edge, $e$, of $\Gamma$ such that the splitting of $G$
which the edge $e$ gives, along the edge group, $E$, is
hyperbolic-hyperbolic with respect to the splitting along $C_3$.
The group $E$  acts hyperbolically on $T_3$. $F$ denotes the fiber
group of $S(C_1,C_2)$. We have the following lemma:
\begin{lem}[Elliptic fiber]
\label{lem,fiber}
Letting $E$ act on  $T_3$,
we obtain a presentation
$$E=\langle t,F|tFt^{-1}=\a(F) \rangle,$$
where $\a \in {\rm Aut}(F)$ or
$$E=L*_F M,$$
where $[L:F]=[M:F]=2$.
\end{lem}

\proof
We first show that $F$ is elliptic on $T_3$.
To argue by contradiction,
assume  that there is $a\in F$ acting hyperbolically on $T_3$. Then the
pairs
$C_1,C_3$ and $C_2,C_3$ are hyperbolic-hyperbolic. Let $F_1$ be the
fiber of the enclosing
group $S(C_1,C_3)$ corresponding to the pair
 $C_1,C_3$ and let $F_2$ be the fiber of $S(C_2,C_3)$.
We claim that there is $w_1\in F_1$ which does not lie
in any conjugate of  $F$. Indeed $F_1$ is contained
in a conjugate of $C_1$. $C_1$ acts on $T_2$ hyperbolically preserving
an axis which is stabilized by $F$.
If $F_1$
contains an element, $w_1$, that acts hyperbolically on $T_2$ then $w_1$
does not lie in any conjugate of $F$. Otherwise $F_1$ fixes an axis and
it
is contained in a conjugate of $F$.  In this case consider the
actions of $C_1$ on $T_2$ and $T_3$. By passing, if necessary (in the
dihedral
action case),
to a subgroup of index 2 we can assume that $C_1$ is generated by
$\langle t,F \rangle $
where $t$ acts hyperbolically on $T_2$.  Similarly $C_1$ is generated by
some
$x$ acting hyperbolically on $T_3$ and a conjugate of $F_1$ which is
contained
in $F$. Since $C_1$ acts hyperbolically on $T_2$ $x$ acts hyperbolically
on
$T_2$ and $x=tf$ where $f\in F$. Since $a$ acts hyperbolically on
$T_3$ we have
$a=x^kf'$ with $f'\in F$.Then $t^{-k}a$ acts elliptically on $T_3$,
therefore
it lies in $F$. But this is a contradiction since $t\notin F$.

Now if $b\in C_3$ either $b\in F_2$ or $b^kw_1^n\in F_2$. This is
because
if $b\notin F_2$,
$w_1,b$ act both
as hyperbolic elements on $T_2$
and they fix the same axis (since $w_1,b \in C_3$ and $C_3$ is slender).
But then $b^k\in S(C_1,C_2)$.
Since the translation length of any hyperbolic element of $C_3$, for its
action
on $T_2$, is a multiple of a fixed number we can pick
the same $k$ for all $b\in C_3$. So one has $C_3^k\subset S(C_1,C_2)$.
Therefore if we consider
the graph of groups corresponding to $S(C_1,C_2)$ and its Bass-Serre
tree
then $C_3$ fixes
a vertex of this tree.
Therefore either it fixes the vertex stabilized by  $S(C_1,C_2)$
or $C_3^k$ is contained in the edge stabilizer of an edge
adjacent to the vertex stabilized
by $S(C_1,C_2)$.
But in the first case we have that $C_3\subset S(C_1,C_2)$ and in the
second
it is impossible that $C_3$ is hyperbolic-hyperbolic
with respect to, say, $C_1$.  We conclude that there is no $a\in F$
acting
hyperbolically on $T_3$. Therefore $F$ is elliptic on $T_3$.

On the other hand the splitting over $C_3$ is hyperbolic-hyperbolic with
respect one of the splittings
used to construct $\Gamma $.
Let's say that it is hyperbolic-hyperbolic with respect to
the splitting over $C_1$. Since $F\subset C_1 $ and $F$ fixes an axis
of $T_3$ a conjugate
of $F$ is contained in $C_3$. Therefore since the splitting over $C_3$ is
hyperbolic-hyperbolic
with respect to the splitting over $E$, $E$ contains a conjugate of
$F$.
Moreover this conjugate of $F$ is  an infinite index subgroup of $E$.
This clearly implies that
$$E=\langle t,F |tFt^{-1}=\a (F) \rangle $$
where $\a $ is an automorphism of $F$ or that
$$E=L*_F M,$$
where $[L:F]=[M:F]=2$. \qed(Lemma \ref{lem,fiber}).

Let $\{e_i\}$be the collection of the edges of $\C$ whose edge
groups, $E_i$, are hyperbolic on $T_3$, and $\{d_j\}, \{D_j \}$
the collections of the rest of the edges and their edge groups.
Let $T_{E_i}$ be the Bass-Serre tree of the splitting of $G$ we
obtain by collapsing all edges of $\C$ but $e_i$. The group $C_3$
is hyperbolic on $T_{E_i}$ by the way we took $e_i$, and the
splitting along $C_3$ is minimal. Consider the diagonal action of
$G$ on $T_3 \times T_{\C}$. In the same way as in the proof of
Proposition \ref{prop,corecomplex},
 we can show that there is
a subcomplex $ \tilde Z$ of $T_3 \times T_{\C}$ which is invariant
by $G$ such that $\tilde Z/G$ is finite. We explain this in
detail: Let $S$ be the enclosing group of $\C $ and let $T_S$ be
the minimal invariant subtree of $T_3$ for the action of $S$. For
each $E_i\subset S$ let $l_i$ be the invariant line for the action
of $E_i$ on $T_3$. Finally for each $D_j$ we pick a vertex $v_j$
on $T_3$ fixed by $D_j$. Let $\tilde e_i$ be a lifting of $e_i$ to
$T_1\times T_2$ with an endpoint on $l_i$ and $\tilde d_j $ a
lifting of $d_j $ on $T_1\times T_2$ with an endpoint on $v_j$.
Let $$Z_1=T_S\cup_i (l_i\times \tilde e_i) \cup_j (\tilde d_j )$$
where the union is over all the $\tilde e_i's, \tilde d_j's$. We
take then $Z$ to be the complex obtained by the translates $GZ_1$.

 As before, we give a description of $Z$
using gluings of graphs. Let $\{k_i\}$ be the vertices of $\Gamma$
other than the one for $S$, and $\{K_i\}$ their vertex groups. Let
$T_3/S={\mathcal S}, T_3/E_i={\mathcal E}_i, T_3/D_i={\mathcal
D}_i, T_3/K_i={\mathcal K}_i$ be the quotient graph of groups.
Since the action of $E_i$ on $T_3$ is hyperbolic and $E_i$ is
slender, there is an invariant line $l_i$ in $T_3$ by $E_i$ and
the core of ${\mathcal E}_i$ is $c_i=l_i/E_i$, which is
topologically a segment if the action of $E_i$ on $T_3$ is
dihedral, or else a circle. A core of ${\mathcal D}_i$ is a vertex
since the action is elliptic. Let's denote a core of a graph of
groups, ${\mathcal A}$, as ${\rm co}({\mathcal A})$. A core
complex, $Z$,
 of the diagonal action of $G$ on $T_{\C} \times T_3$ is
given as follows:
$$Z={\rm co}({\mathcal S})
\bigcup \cup_i{\rm co}({\mathcal K}_i)
\bigcup \cup_i({\rm co}({\mathcal E}_i) \times [0,1]) \bigcup
\cup_i({\rm co}({\mathcal D}_i )
\times [0,1]).$$
Note that each ${\rm co}({\mathcal E}_i) \times [0,1]$
and ${\rm co}({\mathcal D}_i)  \times [0,1]$ is attached to
${\rm co}({\mathcal S}) \bigcup \cup_i{\rm co}({\mathcal K}_i)$
by the graph morphism induced by the homomorphism
of each of $E_i,D_j$
to $S$ and to $K_k$ given in $\Gamma$.

$Z$ is a finite complex, and the fundamental group in the sense of
complexes of groups (let's call such fundamental group {\it
H-fundamental group} in this proof) is $G$. Let ${\mathcal
C}_3=T_{\Gamma}/C_3$. Since $C_3$ is slender, and the action of
$C_3$ on $T_{\Gamma}$ is hyperbolic, there is an invariant line
$l$ in $T_{\Gamma}$. Since the splittings of $G$ along $E_i,C_3$
are minimal, we can conclude, as in Prop.
\ref{prop,enclosinggroup}, that $\cup_i({\rm co}({\mathcal E}_i)
\times [0,1]) ={\rm co}({\mathcal C}_3)\times [0,1]$
 in $(T_{\Gamma} \times T_3)/G$.
Although $l/C_3$ is embedded in $Z$, which locally separates $Z$,
the splitting of $G$ along $C_3$ which we obtain by cutting $Z$
along $l/C_3$ might be different from the original splitting along
$C_3$.

Consider the following
subcomplex, $W$, of $Z$,
$$W={\rm co}({\mathcal S}) \bigcup
\cup_i({\rm co}({\mathcal E}_i) \times [0,1]) \bigcup \cup_i({\rm
co}({\mathcal D}_i) \times [0,1]).$$ Let $\{p_j\}$ be the set of
vertices in $W$ which are not contained in ${\rm co}({\mathcal
S})$. Let $m_j$ be the link of $p_j$ in $W$, which we denote by
$Lk(p_j,W)$. Since each ${\rm co}({\mathcal D}_i)$ is a point, if
$p_j$ is in $\cup_i({\rm co}({\mathcal D}_i) \times [0,1]) $, then
$m_j$ is a point, whose fundamental group (in the sense of graph
of groups) is one of the $D_i$'s (the group corresponding to the
edge which contains $p_j$). The point $m_j$ locally separates $W$,
and also $Z$. If the vertex $p_j$ is in $\cup_i({\rm co}({\mathcal
E}_i) \times [0,1])$, then the link $m_j$ is the finite union of
circles and segments, such that each of them locally separates
$W$, and also $Z$.

If we cut $Z$ along the union of those links $\cup_j Lk(p_j,W)$,
we obtain a graph decomposition of $G$ by Lemma
\ref{lem,vankampen} such that edge groups are the image in $G$ of
the H-fundamental groups of connected components of $\cup_j
Lk(p_j,W)$.
Let $V$ be the connected component of $W - \cup_j Lk(p_j,W)$ which
contains ${\rm co}({\mathcal S})$. The image in $G$ of the
H-fundamental group of $V$ contains $S=S(C_1,C_2)$. Let's denote
it by $S'$. We claim that $S'$ is an extension of the fundamental
group of some $2$-orbifold, $\Sigma'$, by $F$, the fiber group of
$S$ such that $\Sigma \subset \Sigma'$. To see it, let $U=\cup_i
({\rm co}({\mathcal E}_i)\times [0,1])$, which is a squared
surface possibly with some vertices identified. Note $U \subset
W$. Let $\{q_i\}$ be the vertices in $U \cap {\rm co}({\mathcal
S})$. Define $l_i=Lk(q_i,U)$ for each $i$. Note that each $p_i \in
U$ and also $m_i \subset U$. If we cut $U$ along $\cup_i l_i$ and
$\cup_i m_i$, we obtain a graph decomposition along slender
groups, which are the image (in $G$) of the H-fundamental group of
$l_i$'s and $m_i$'s. Let $U' \subset U$ be the connected component
of $U-(\cup_i l_i \bigcup \cup_i m_i)$ which does not contain any
of $p_i,q_i$, i.e., the vertices of $U$. We know that $U'$ is a
surface with boundary. Also $U' \subset V$. Cutting $V$ along
$\cup_i l_i$, where $U'$ is one of the connected component after
the cutting, we obtain a graph decomposition of $S'$ along the
slender groups corresponding to $l_i$'s. The vertex group, $S_0$,
corresponding to $U'$ is an extension of the fundamental group of
some $2$-orbifold, $\Sigma_0$, by $F$ by our construction.
$\Sigma_0$ is obtained from the $2$-manifold $U'$ attaching a
disks or a half disk with cone points appropriately each time if
the fundamental group of $m_i$ does not inject to $G$ (cf. we did
the same thing when we constructed $\Sigma$ for $S$ previously). A
vertex group other than $S_0$ is not only a subgroup of $S$, but
also it corresponds to a suborbifold in $\Sigma$, the $2$-orbifold
for $S$. To see it, consider a small neighborhood, $\overline{{\rm
co}({\mathcal S})}$, of ${\rm co}({\mathcal S})$ in $W$. To be
concrete, for example, we take a barycentric subdivision of $W$
and collect all cells which intersect ${\rm co}({\mathcal S})$.
The H-fundamental group of $\overline{{\rm co}({\mathcal S})}$ is
$S$. One can consider that ${\rm co}({\mathcal S})$ is a
deformation retract of $\overline{{\rm co}({\mathcal S})}$. We may
assume that each $l_i$ is in $\overline{{\rm co}({\mathcal S})}$.
Cutting $\overline{{\rm co}({\mathcal S})}$ along $\cup_i l_i$, we
obtain a graph decomposition of $S$ along slender groups.
 This decomposition is realized by cutting $\Sigma$
along simple closed curves and segments. (Consider
the quotients by $F$ of the H-fundamental
groups of $\overline{{\rm co}({\mathcal S})}$ and
$\l_i$'s and obtain a decomposition
of the orbifold fundamental group of $\Sigma$ along slender groups,
and reduce the argument to surface topology.
Note that all maximal peripheral subgroups
of $S$ are elliptic with respect to the graph decomposition,
cf. Rem \ref{rem,maxperipheral}, so that $\Sigma$ does
not have any free boundary points in the decomposition, and
the H-fundamental group of $l_i$ injects in $G$).
Let $S_i$ be the H-fundamental group
of the connected component of $\overline{{\rm co}({\mathcal S})} - \cup_i l_i$
which contains $q_i$.
Note that this is identical to the connected component
of $V-\cup_i l_i$ which contains $q_i$.
Let $\Sigma_i \subset \Sigma$ be the sub-orbifold
such that $S_i$ is the H-fundamental group
of $\Sigma_i$ (each $\Sigma_i$ is a connected component
of $\Sigma$ after the cutting we obtained in the above).
Then, $S_i$ is an extension of the orbifold fundamental
group of $\Sigma_i$ by $F$.
Since the graph decomposition of $S'$ we obtained
by cutting $V$ along $\cup_i l_i$ has vertex groups $S_i$'s
(corresponding to $q_i$'s) and
$S_0$, with edge groups corresponding to the H-fundamental
groups of $l_i$'s, and
$S_0$ is also an extension
of the orbifold fundamental group of $\Sigma_0$ by $F$, we
conclude that $S'$ is an extension of the orbifold fundamental
group of a $2$-orbifold, $\Sigma'$, by $F$ such that
$\Sigma'$ is the union of $\Sigma_i$'s and $\Sigma_0$
pasted along $l_i$'s. By construction, $\Sigma \subset \Sigma'$.

Let $\Gamma'$ be the graph decomposition of $G$ obtained by cutting $Z$ along
$\cup_j Lk(p_j,W)$, with a vertex group $S'$.
We first show that $\Gamma'$ satisfies the
properties 1,2,3 of Def \ref{defn,enclosing} for the three splittings
(cf. Prop \ref{prop,enclosinggroup}).
Then we apply Prop \ref{prop,modification} to $\Gamma'$ and obtain a refinement,
$\Gamma''$, such that each edge of $\Gamma''$ gives a minimal splitting of $G$.
We will show that $\Gamma''$ satisfies all the properties
of Def \ref{defn,enclosing} so that it is an enclosing graph decomposition
for the three splittings, with an enclosing vertex group $S'$.
The argument is similar to Prop
\ref{prop,enclosingdecomposition}.

The group $S'$ has the property 3 by the construction.
Regarding the property 1 of $\Gamma'$, it is clear that $S'$ contains some
conjugates of $C_1,C_2,C_3$.

To verify the property 2(rigidity) of $S'$, let $\Lambda$ be a graph decomposition
of $G$ such that the splitting of $G$ which any edge of $\Lambda$ gives
is elliptic-elliptic with respect to any of the three splittings.
We argue in the same way as
in the proof of Prop \ref{prop,enclosinggroup}
to show that  $S'$ is elliptic on the Bass-Serre tree
of $\Lambda$, $T_{\Lambda}$.
Let $c$ be either an essential simple closed curve on $\Sigma'$
or an essential embedded segment on $(\Sigma', (\partial \Sigma';2))$.
Let $C<S'$ be the fundamental group for $c$.
Cutting $\Sigma'$ along $c$, we obtain a splitting of $G$ along $C$.
This  splitting is minimal by Prop \ref{prop,dual-minimal}.
Let $T_C$ be its Bass-Serre tree.
For our purpose, by Cor \ref{cor,elliptic},
it suffices for us to show that the group $C$
is elliptic on $T_{\Lambda}$ to conclude that so is $S'$.
Let $e$ be an edge of $\Lambda$, which gives a splitting
of $G$ along its edge group, $E$. To conclude $C$ is elliptic on $T_{\Lambda}$,
we will show that $C$ is elliptic with respect to the splitting
along $E$.
Since the group $E$ is elliptic with respect
to the (original) splittings of $G$ along $C_1,C_2$,
it is elliptic on $T_{C_1} \times T_{C_2}$, i.e., $E$ fixes a vertex.
Therefore $E$ is elliptic on $T_{\Gamma}$, which is the Bass-Serre tree
of the enclosing graph decomposition, $\Gamma$,
 we constructed for the splittings
along $C_1,C_2$.  Moreover, we know that $E$ is not in a conjugate of $S$ (cf.
the proof of Prop \ref{prop,enclosinggroup}).
Since the group $E$ is elliptic on $T_{C_3}$ as well by our assumption,
it fixes a vertex when it acts on $T_{\Gamma} \times T_{C_3}$.
It follows that $E$ is elliptic on $T_{\Gamma'}$, the Bass-Serre tree
for $\Gamma'$,  by the way we constructed it. Therefore
$E$ is in a conjugate of a vertex group of $\Gamma'$, which
is not $S'$. This implies that the group $E$ is elliptic on $T_C$.
Since the splitting along $C$ is minimal, we find that
the pair of splittings along $C$ and $E$ is elliptic-elliptic.
But, the edge $e$ was an arbitrary edge of $\Lambda$, so that
the group $C$ is elliptic on $T_{\Lambda}$.
We showed the property 2 for $S'$.

So far, we have shown that the graph decomposition $\Gamma'$ with
a vertex group $S'$ satisfies the properties 2,3 and a part of the
property 1.
As we obtain Prop \ref{prop,enclosinggroup} from Prop
\ref{prop,enclosingdecomposition} for a pair
of splittings, we apply Prop \ref{prop,modification}
to $\Gamma'$ and obtain a graph decomposition $\Gamma''$
such that each edge gives a minimal splitting.
We now claim that $\Gamma''$
has $S'$ as an enclosing vertex group and satisfies all
the properties to be an enclosing decomposition for the
three splittings.
The argument is same as when we show Prop \ref{prop,enclosingdecomposition} from
Prop \ref{prop,enclosinggroup},
so we omit some details.
By construction, $\Gamma''$ has the property 4.
$S'$ is a vertex group of $\Gamma''$ because the edge
groups of $\Gamma''$ are in edge groups of $\Gamma'$ and
the rigidity of $S'$.
Therefore, $\Gamma''$ with $S'$ satisfies the properties 2, 3,
and the property 1 except for the last item, which
we did not verify for $\Gamma'$.

To verify the rest of the property 1 for $\Gamma''$,
let $e$ be an edge of $\Gamma''$ with the edge group, $E$.
Let $T_{C_i}$ be the Bass-Serre tree of the (original) splitting
of $G$ along $C_i$, $i=1,2,3$.
We want to show that the edge group $E$ is elliptic on all
$T_{C_i}$.
The splitting of $G$ along $E$ which the edge $e$
gives is minimal.
Let $T_E$ be the Bass-Serre tree of this splitting.
By the property 2 (rigidity) of $S'$,
$S'$  is elliptic on $T_E$, so that the subgroups
$C_i$ are elliptic as well. It follows that the group $E$
is elliptic on all $T_{C_i}$ because the splitting along $E$
is minimal. This is what we want.
We showed all the properties for $\Gamma''$ with $S'$, so that
the proof of Prop \ref{prop,enclosing} is complete.
\qed
(Prop \ref{prop,enclosing})

\subsection{Maximal enclosing decompositions}
Following the previous subsection, we produce an enclosing graph
decomposition of a set, $I$,  of hyperbolic-hyperbolic minimal
splittings of $G$ along slender subgroups. We put an order to the
elements in $I$ such that if $I_i$ denotes the set of the first
$i$ elements, then each $I_i$ is a set of hyperbolic-hyperbolic
splittings. Then we produce a sequence of graph decompositions,
$\Gamma_i$, of $G$ such that $\Gamma_i$ is an enclosing graph
decomposition for $I_i$ with enclosing vertex group $S_i$. We
already explained how to construct $\Gamma_2$, then $\Gamma_3$
using it.In the same way as we produce $\Gamma_3$ from $\Gamma_2$
fromthe splitting along $C_3$, we produce $\Gamma_{i+1}$ from
$\Gamma_i$.  Note that $\Gamma_{i+1}$ is identical to $\Gamma_i$
if the edge group, $C_{i+1}$, of the $(i+1)$-th splitting is
contained in a conjugate of $S_i$.

Although $\Gamma_i$ is an infinite sequence in general, there
exists a number $N$ such that $\Gamma_i$ is identical if $i \ge N$
by the following result. We recall that a graph of groups, whose
fundamental group is $G$, is {\it reduced} if its  Bass-Serre tree
does not contain any proper subtree which is $G$-invariant, and
the vertex group of any vertex of the graph of valence $2$
properly contains the edge groups of the associated edges.

\begin{thm}[Bestvina-Feighn accessibility \cite{BF}]
\label{thm,bestvinafeighn} Let $G$ be a finitely presented group.
Then there exists a number $\gamma(G)$ such that if $\Gamma$ is a
reduced graph of groups with fundamental group isomorphic to $G$,
and small edge groups, then the number of vertices of $\Gamma$ is
at most $\gamma(G)$.
\end{thm}

Let $\Sigma_i$ be the $2$-orbifold for the enclosing vertex group
$S_i$. Then, $\Sigma_i \subset \Sigma_{i+1}$ as $2$-orbifolds. Any
system, ${\mathcal F}$, of disjoint essential simple closed curves
on $\Sigma_i$ and essential segments on $(\Sigma_i,(\partial
\Sigma_i;2))$ such that any two of them are not homotopic to each
other gives a reduced graph decomposition of not only $S_i$ but
also $G$. Because the number of the connected components of
$\Sigma_i \backslash {\mathcal F}$ are bounded by $\gamma(G)$ by
Theorem \ref{thm,bestvinafeighn}, there exists $N$ such that
$\Sigma_i$ is constant if $i \ge N$. This implies, by the way we
constructed $\{ \Gamma_i \}$, $\Gamma_i$ is also constant if $i
\ge N$. $\Gamma_N$ is an enclosing graph decomposition for $I$. We
have shown the following.

\begin{prop} \label{prop ,enclosehyperbolic}
Let $G$ be a finitely presented group. Let $I$ be a set of
hyperbolic-hyperbolic minimal splittings of $G$ along slender
subgroups. Then,
 an enclosing graph decomposition of $G$, $\Gamma_I$, exists for $I$.
\end{prop}

If a set of hyperbolic-hyperbolic minimal
splitting along slender groups, $I$, is maximal, we call $\Gamma_I$ maximal.
A maximal enclosing graph decomposition has the following property.

\begin{lem}[Maximal enclosing graph decomposition]
\label{lem,maximalenclosing} Let $G$ be a finitely generated
group. Let $\Gamma$ be a maximal enclosing decomposition of $G$.
Suppose $A*_C B, A*_C$ is a minimal splitting of $G$ along a
slender subgroup $C$. Then the splitting along $C$ is
elliptic-elliptic with respect to the splitting of $G$ which each
edge of $\Gamma$ gives. In particular, the group $C$ is elliptic
on $T_{\Gamma}$, the Bass-Serre tree for $\Gamma$.
\end{lem}

\begin{rem}
Although the existence of maximal enclosing group
is guaranteed only for a finitely presented group,
the lemma is true if a maximal enclosing group
exists for a finitely generated group.
\end{rem}

\proof Suppose not. Then there exists an edge, $e$, of $\Gamma$
with edge group, $E$, such that the minimal splitting of $G$ the
edge $e$ gives is hyperbolic-hyperbolic  with respect to the
splitting along $C$. Suppose $\Gamma$ is an enclosing
decomposition for a maximal set $I$. Since $E<S$, $S$ is
hyperbolic on $T_C$, the Bass-Serre tree for the splitting along
$C$. Then by the rigidity (Def \ref{defn,enclosing}) of $S$, there
is a splitting in $I$ which is hyperbolic-hyperbolic with respect
to the splitting along $C$.
 Let $I'$ be the union
of $I$ and the splitting along $C$, which is a set
of hyperbolic-hyperbolic splittings. If we produce
an enclosing decomposition for $I'$ using $\Gamma$ and the decomposition
along $C$, we obtain a graph decomposition with a  different enclosing
vertex group
(namely, the $2$-orbifold is larger) from $\Gamma$, because the group $C$ is
hyperbolic with respect to the splitting along $E$,
 which is impossible since $\Gamma$ is maximal.
The last claim is clear from Bass-Serre theory. \qed



\begin{prop}[Rigidity of maximal enclosing group]
\label{prop,rigidenclosing} Let $G$ be a finitely generated group.
Let $\Gamma, \Gamma'$ be maximal enclosing decompositions of $G$
with enclosing groups $S, S'$. Then the group $S$ is elliptic on
$T_{\Gamma'}$, the Bass-Serre tree of $\Gamma'$, so that $S$ is a
subgroup of a conjugate of a vertex group of $\Gamma'$. If $S'$ is
a subgroup of a conjugate of $S$, then it is a conjugate of $S$.
\end{prop}

\proof Let $I,I'$ be maximal sets of hyperbolic-hyperbolic minimal
splittings of $G$ for $\Gamma,\Gamma'$. We may assume $I \not=I'$.
Since they are maximal, if they have a common splitting, then
$I=I'$, so that $I \cap I' =\emptyset$. Also, a pair consisting of
any splitting in $I$ and any splitting in $I'$ is
elliptic-elliptic. Moreover, there is no (minimal) splitting of
$G$ along a slender subgroup which is hyperbolic-hyperbolic with
respect to some splittings in both of $I,I'$, because, then such
splitting and $I\cup I'$ would violate the maximality of $I$.

Let $\Sigma$ be the $2$-orbifold for the enclosing group $S$.
Let $d$ be a simple closed curve on $\Sigma$ or a segment on 
$(\Sigma, (\partial \Sigma;2))$
which is essential, then
cutting $\Sigma$ along $d$, we obtain a splitting of $G$
along the group, $D$, which is the fundamental group of $d$.
$D$ is slender, and the splitting of $G$ along $D$ is minimal
since it is hyperbolic-hyperbolic with respect to one
of the minimal splittings in $I$ (Prop \ref{prop,dual-minimal}).
Therefore, the splitting along $D$ is elliptic-elliptic to
all splittings in $I'$.

By Cor \ref{cor,elliptic}, it suffices for us to show that the
group $D$ is elliptic on $T_{\Gamma'}$, the Bass-Serre tree of
$\Gamma'$ to show that $S$ is elliptic on it. Let $e$ be an edge
of $\Gamma'$, with edge group $E$. Collapsing all edges of
$\Gamma'$ except $e$, we obtain a splitting of $G$ along $E$,
which is minimal. Let $T_E$ be the Bass-Serre tree of this
splitting. Then, it is enough for us to show that the group $D$ is
elliptic on $T_E$. Since the splitting along $E$ is minimal, it
suffices to show that the group $E$ is elliptic on $T_D$, the
Bass-Serre tree for the splitting along $D$. Since $E<S'$, it is
enough if we show that $S'$ is elliptic on $T_D$. By the property
2 (rigidity) of $S'$, it suffices to show that the splitting along
$D$ is elliptic-elliptic with respect to all splittings in $I'$,
which we already know. We have shown that $S$ is elliptic on
$T_{\Gamma}$.

By the same argument, $S'$ is elliptic on $T_{\Gamma}$. Suppose
$S$ is in a conjugate of $S'$, i.e., $S < gS'g^{-1}, g \in G$.
Then $S'$ is also in a conjugate of $S$, since, otherwise, $S'$ is
in a conjugate of a vertex group of $\Gamma$ which is not $S$.
Then this vertex group contain a conjugate of $S$, which is
impossible since all edge of $\Gamma$ which is adjacent to the
vertex whose vertex group is a conjugate of $S$ has an edge group
which is a proper subgroup of the conjugate of $S$. Suppose $S' <
hSh^{-1}, h \in G$. Therefore, $S < ghS(gh)^{-1}$. This implies
that $gh \in S$, and $S=ghS(gh)^{-1}$, so that $S = gS'g^{-1}$.
\qed

\subsection{JSJ-decomposition for hyperbolic-hyperbolic minimal splittings}
Using maximal enclosing graph decompositions, we produce
a graph decomposition of $G$, $\Lambda$,
which "contains" all maximal enclosing groups.
$\Lambda$ will deal with all minimal splittings
of $G$ along slender subgroups which are hyperbolic-hyperbolic
with respect to some (minimal) splittings along slender subgroups.


Consider all maximal enclosing decompositions, $\C_i$, of $G$ with
enclosing groups, $S_i$. Let $T_i$ be the Bass-Serre tree of
$\C_i$. We construct a sequence of refinements $\{\Lambda_i\}$
such that $\Gamma_1=\Lambda_1$. We then show that after a finite
step, the graph decompositions stay the same. We denote the
decomposition obtained after this step by $\Lambda$.

We put $\Lambda_1=\Gamma_1$. We consider now $\Gamma_2$. By Prop
\ref{prop,rigidenclosing}, $S_2$ is elliptic on $T_1$. If $S_2$ is
a subgroup of a conjugate of $S_1$, then we do nothing and put
$\Lambda_2=\Gamma_1$. If $S_2$ is conjugate into a vertex group,
$A$, of $\C_1$ which is not $S_1$, then we let $A$ act on $T_2$
and obtain a refinement, $\C_2'$, of  $\C_1$. Namely, let $\A$ be
the graph decomposition of $A$ we get. We substitute $\A$ to the
vertex, $a$, for $A$ in $\Gamma_1$ (see Def \ref{defn,refinement} and
the following remarks). We can do this since each edge group, $E$,
of $\C_1$ is elliptic on $T_2$ (because $E$ is a subgroup of
$S_1$, which is elliptic on $T_2$), so that $E$ is a subgroup of a
conjugate of a vertex group of $\A$.  Note that all edge groups of
$\Gamma_2'$ are slender since they are subgroups of conjugates of
edge groups of $\Gamma_1, \Gamma_2$. $\C_2'$ has conjugates of
$S_1,S_2$ as vertex groups. Also they are peripheral subgroups of
either $S_1$ or $S_2$.


Each edge, $e$, of $\Gamma'_2$ gives a minimal splitting of $G$
along its edge group, $E$, which is a subgroup of a conjugate of
the edge group, $E'$, of an edge, $e'$, of $\Gamma_i, 
(i=1 \, {\rm or }\, 2)$. 
We show this by contradiction: suppose that the edge $e$
gives a non-minimal splitting, which is hyperbolic-elliptic with
respect to a splitting $G=P*_D Q, ({\rm or } P*_D)$ such that $D$
is slender. Let $T_D$ be its Bass-Serre tree. Since the group $E$
is hyperbolic on $T_D$, so is $E'$. Because the splitting of $G$
along $E'$, the one $e'$ gives, is minimal, it is
hyperbolic-hyperbolic with respect to $G=P*_D Q, ({\rm or }
P*_D)$. By Prop \ref{prop,dual-minimal}, the splitting along $D$
is minimal. On the other hand,  by Lemma
\ref{lem,maximalenclosing}, the group $D$ is elliptic on
$T_{\Gamma_i}$, the Bass-Serre tree of $\Gamma_i$ since it is
maximal. It follows that the group $D$ is elliptic on $T_{E'}$,
the Bass-Serre tree of the splitting along $E'$, a contradiction.


We collapse all edges of this decomposition which are not adjacent
to the vertices with vertex groups $S_1,S_2$. If the resulting
graph decomposition is not reduced (cf. Theorem
\ref{thm,bestvinafeighn}) at some vertex, then we collapse one of
the associated two edges, appropriately, to make it reduced. Note
that it is reduced at a vertex whose vertex group is an enclosing
group since all edge groups are proper subgroups at the vertex of
an enclosing group. We denote the resulting reduced graph
decomposition by $\Lambda_2$. We remark that $\Lambda_2$ is a
refinement of $\Lambda_1$.

By our construction, all edge groups of $\Lambda_2$ are
conjugates of edge groups of $\Gamma_1,\Gamma_2$,
and each edge of $\Lambda_2$  is connected to the vertex of an enclosing group.
This is not obvious, we recall that we only know that edge groups
of $\C_1 $ are elliptic on $T_1$. When we substitute ${\mathcal
A}$ for $A$, some of these edge groups are connected to (a
conjugate of) $S_2$.
We have to show that these edge groups are peripheral in
$gS_2g^{-1}$. To see it, let $E$ be an edge group, and suppose
that $E<S_2$, where we assume that $g=1$ for notational
simplicity. (In general, just take conjugates by $g$ appropriately
in the following argument). We will show $E$ is peripheral in
$S_2$. Note that this is the only essential case since $E$ can be
only peripheral in $S_1$ because we don't do anything around the
vertex for $S_1$ when we construct $\Lambda_2$. Also we may assume
$E<S_1$. Let $\Sigma_2$ be the $2$-orbifold for $S_2$, and $d$ an
essential simple closed curve/segment on it with the group $D$
represented by $d$. Cutting $\Sigma_2$ along $d$, we obtain a
splitting of $G$ along $D$, with Bass-Serre tree $T_D$. It
suffices to show that the group $E$ is elliptic with respect to
$T_D$ to conclude that $E$ is peripheral in $S_2$. Suppose not.
Then, the splittings of $G$ along $E$ and $D$ are
hyperbolic-hyperbolic since both of them are minimal. Then $S_1$
is hyperbolic on $T_D$ since $E<S_1$. Let $\Sigma_1$ be the
$2$-orbifold for $S_1$. It follows from Cor \ref{cor,elliptic}
 that there exists an essential simple closed curve or
a segment, $d'$,  on $\Sigma_1$ such that cutting $\Sigma_1$ along
$d'$ gives a splitting of $G$ along the group for $d'$, $D'$, such
that the splittings along $D$ and $D'$ are hyperbolic-hyperbolic.
Then, the set $I_1 \cup I_2$ with the two splittings along $D,D'$
is a set of hyperbolic-hyperbolic, which is impossible
because the enclosing group for this set must be strictly bigger
than $S_1$, and $S_2$ as well, which is impossible since they
are maximal. We have
show that all edge groups of $\Lambda_2$ are peripheral subgroups
of enclosing vertex groups.

We continue similarly and obtain a sequence of reduced graph
decompositions of $G$; $\Lambda_1, \Lambda_2, \Lambda_3, \cdots$.
Namely, we first show that the enclosing group $S_3$ is elliptic
with respect to $\Lambda_2$, using the maximality of $S_i$ and
rigidity. If $S_3$ is a subgroup of a conjugate of $S_1$ or $S_2$,
then $\Lambda_3$ is $\Lambda_2$. Otherwise, there exists a vertex
group, $A_2$, of $\Lambda_2$ which is different from $S_1,S_2$ and
contains a conjugate of $S_3$. We let $A_2$ act on the Bass-Serre
tree of $\Gamma_3$ and obtain a graph decomposition, $\A_2$, of
$A_2$. We then substitute $\A_2$ to the vertex for $A_2$ in
$\Lambda_2$, which is $\Gamma'_3$. We show that all edges of
$\Gamma'_3$ give minimal splittings of $G$ along slender
subgroups. We then collapse all edges which are not adjacent to
the vertices with the vertex groups conjugates of $S_1,S_2,S_3$,
and also collapse edges appropriately at non-reduced vertices, to
obtain a reduced graph decomposition, $\Lambda_3$. The edge groups
of $\Lambda_3$ which are adjacent to some conjugates of $S_i$ are
the conjugates of peripheral subgroups of $S_i$. In this way, we
obtain $\Lambda_{n+1}$ from $\Lambda_n$ using $\Gamma_{n+1}$. This
is a sequence of refinements.

We claim that there exists a number $N$ such that if $n \ge N$
then $\Lambda_{n+1}=\Lambda_n$. Indeed, if not, then
the number of vertices in $\Lambda_n$ whose
vertex groups are enclosing groups $S_i$ tends to infinity
as $n$ goes to infinity. This is impossible since
the number of the vertices of $\Lambda_n$ is at most $\gamma(G)$ by
Theorem \ref{thm,bestvinafeighn}.
Note that slender groups are small.
Let's denote $\Lambda_N$ by $\Lambda$, and
state some of the properties we have shown as follows.

\begin{prop}[JSJ-decomposition for hyp-hyp
minimal splittings along slender groups]
\label{prop,jsj.hyp-hyp.minimal} Let $G$ be a finitely presented
group. Then there exists a reduced graph decomposition, $\Lambda$,
with the following properties:
\begin{enumerate}
\item all edge groups are slender. \item Each edge of $\Lambda$
gives a minimal splitting of $G$ along a slender group. This
splitting is elliptic-elliptic with respect to any minimal
splitting of $G$ along a slender subgroup.
 \item Each maximal
enclosing group of $G$ is a conjugate of some vertex group of
$\Lambda$, which we call a (maximal) enclosing vertex group. The
edge group of any edge adjacent to the vertex of a maximal
enclosing vertex group is a peripheral subgroup of the enclosing
group.
 \item Each edge of $\Lambda$ is adjacent to some vertex
group whose vertex group is a maximal enclosing group.
\item Let
$G=A*_C B$ or $A*_C$ be a minimal splitting of $G$ along a slender
subgroup $C$, and $T_C$ its Bass-Serre tree.
\begin{enumerate}
\item
If  it is hyperbolic-hyperbolic
with respect to a minimal splitting of $G$ along a slender
subgroup, then
\begin{enumerate}
\item\label{item,hyperbolic.enclosing} a conjugate of $C$ is a
subgroup of a unique enclosing vertex group, $S$, of $\Lambda$.
$S$ is also the only one among enclosing vertex groups which is
hyperbolic on $T_C$. There exists a base $2$-orbifold, $\Sigma$,
for $S$ and an essential simple closed curve or a segment on
$\Sigma$ whose fundamental group (in the sense of complex of
groups) is a conjugate of $C$.
\item Moreover, if $G$ does not
split along a group which is a subgroup of $C$ of infinite index,
then all non-enclosing vertex groups of $\Lambda$ are elliptic on
$T_C$.
\end{enumerate}
\item
If it is elliptic-elliptic with respect to any
 minimal splitting of $G$ along a slender
subgroup, then all vertex groups of $\Lambda$ which are
maximal enclosing groups are elliptic on $T_C$.
\end{enumerate}

\end{enumerate}
\end{prop}

\proof By the previous discussion we know that properties 1,3,4,
and a part of property 2 hold. Let's show the rest of the property
2. To argue by contradiction, suppose that the edge, $e$, of
$\Lambda$ gives a minimal splitting along the edge group, $E$,
 which is hyperbolic-hyperbolic with respect
to a minimal splitting of $G$ along a slender subgroup, $C$. But
then $C$ would be contained in an enclosing vertex group of some
graph decomposition $\Gamma_i$ and it would not be a peripheral
group in $\Gamma_i$, a contradiction

We show now (5-a). There is a maximal enclosing group which
contains a conjugate of $C$, such that $C$ is the fundamental
group of an essential simple closed curve or a segment of the
$2$-orbifold for the enclosing group. This is because we start
with the splitting along $C$ to construct the enclosing group. By
the construction of $\Lambda$, this enclosing group is a conjugate
of some vertex group, $S$, of $\Lambda$. To argue by
contradiction, suppose there is another enclosing vertex group
which contains a conjugate of $C$. Then, by Bass-Serre theory,
there must be an edge associated to each of those two vertices
whose edge group contains a conjugate of $C$. But the edge and its
edge group has the property 2, which contradicts the assumption on
the splitting along $C$ that it is hyperbolic-hyperbolic. One can
show that all enclosing vertex groups of $\Lambda$ except $S$ are
hyperbolic on $T_C$ using Cor \ref{cor,elliptic}, and we omit
details since similar arguments appeared repeatedly.

To show the last claim, suppose that there is a non-enclosing
vertex group, $V$, of $\Lambda$ which is hyperbolic on $T_C$.
Letting $V$ act on $T_C$, we obtain a graph decomposition
of $V$, which we can substitute for $V$ in $\Lambda$.
All edge groups of this graph decomposition are conjugates
of subgroups of $C$, which have to be of finite index by our
additional assumption.
Since a conjugate of $C$ is contained in $S$, which is different from $V$,
by Bass-Serre theory, there must be an edge in $\Lambda$ adjacent to the vertex
for $S$ whose edge group, $E$, is a conjugate of a subgroup of $C$ of
finite index. But the edge and its edge group $E$, so that the group $C$ as well,
satisfies the property 2, which contradicts the assumption on $C$
that it is hyperbolic-hyperbolic.

We show (5-b). Let $S$ be a maximal enclosing vertex group of
$\Lambda$. Let $\Gamma$ be a maximal enclosing decomposition which
has $S$ as the maximal enclosing group. Let $\Sigma$ be the
$2$-orbifold for $S$, and $c$ an essential simple closed curve or
an essential segment on $\Sigma$. Cutting $\Sigma$ along $d$, we
obtain a splitting of $G$ along the slender group, $D$,  which
corresponds to $d$. This splitting is minimal by Prop
\ref{prop,dual-minimal}. By the assumption on the splitting along
$C$, the group $D$ is elliptic on $T_C$. It follows from Cor
\ref{cor,elliptic} that $S$ is elliptic on $T_C$. \qed

\subsection{Elliptic-Elliptic splittings}
Let $G=A_n*_{C_n} B_n$ (or $A_n*_{C_n}$) be all minimal splittings
of $G$ along slender groups, $C_n$, which are elliptic-elliptic
with respect to any minimal splitting of $G$ along slender
subgroups. To deal with them as well, we refine $\Lambda$ which we
obtained in Prop \ref{prop,jsj.hyp-hyp.minimal}. As we constructed
a sequence of refinements $\{\Lambda_n\}$ to obtain $\Lambda$, we
construct a sequence, $\{\Delta_n\}$, of refinements of $\Lambda$
using the sequence of splittings along $C_n$. Then we show that
after a finite step the sequence stabilizes in some sense, again
by Theorem \ref{thm,bestvinafeighn}, and obtain the desired graph
decomposition of $G$, $\Delta$.

We explain how to refine $\Lambda$ in the first step. Let $G=A*_C
B$ or $A*_C$ be a minimal splitting along a slender group $C$
which is elliptic-elliptic with respect to any minimal splitting
of $G$ along a slender group. By Prop
\ref{prop,jsj.hyp-hyp.minimal}, all enclosing vertex groups and
all edge groups of $\Lambda$ are elliptic on $T_C$, the Bass-Serre
tree of the splitting along $C$. Let $U$ be a vertex group of
$\Lambda$ which is not an enclosing vertex group. Letting $U$ act
on $T_C$, we obtain a graph decomposition of $U$, ${\mathcal U}$,
which may be a trivial decomposition. Substituting ${\mathcal U}$
for the vertex for $U$ in $\Lambda$, which one can do since all
edge groups of $\Lambda$ are elliptic on $T_C$, we obtain a
refinement of $\Lambda$. We do this to all non-enclosing vertex
groups of $\Lambda$. Then we apply Prop \ref{prop,modification} to
this graph decomposition, and obtain a further refinement of
$\Lambda$, which we denote $\Delta_1$, such that each edge of
$\Delta_1$ gives a minimal splitting of $G$ along a slender group.
By construction, all vertex groups of $\Delta_1$ are elliptic on
$T_C$. Although when we apply Prop \ref{prop,modification}, a
vertex group may become smaller, all enclosing vertex groups of
$\Lambda$ stay as vertex groups in $\Delta_1$. We see this by an
argument  similar to the one in the proof of Prop
\ref{prop,jsj.hyp-hyp.minimal}. We omit details, but just remark
that all edge groups of $\Delta_1$ are edge groups of $\Lambda$ or
subgroups of conjugates of $C$.

Note that $\Delta_1$ might not be reduced. This may cause a
problem when we want to apply Theorem \ref{thm,bestvinafeighn}
later. To handle this problem, if there is a vertex of $\Delta_1$
of valence two such that one of the two edge group is same as the
vertex group, we collapse that edge. Note that the other edge
group is properly contained in the vertex group in our case. We do
this to all such vertices of $\Delta_1$ at one time, and obtain a
reduced decomposition, which we keep denoting $\Delta_1$. In
general, we can obtain a reduced decomposition from a non-reduced
decomposition in this way. We call the inverse of this operation
an {\it elementary unfolding}. By definition, a composition of
elementary unfoldings is an elementary unfolding. If we obtain a
graph decomposition, $\Gamma'$, by an elementary folding from
$\Gamma$, we may say $\Gamma'$ is an elementary unfolding of
$\Gamma$.

\begin{example}[Elementary unfolding]
\label{example,unfolding}
Let $\C$ be a graph decomposition of $G$
which is $G=A*_P B*_Q C$ and suppose $B$ has a graph decomposition
${\mathcal B}$ which is $P*_{P'} B' *_{Q'} Q$.
Then one can substitute ${\mathcal B}$ to the vertex of $B$ in $\C$
and obtains a new graph decomposition $\C'$, which is an
elementary unfolding.
One can substitute ${\mathcal B}$ to $\C$ because
each edge group adjacent to the vertex for $B$ in $\C$
($P,Q$ in this case) is a subgroup of some vertex group of  ${\mathcal B}$.
If we refine $\C$ using ${\mathcal B}$, we obtain $\C'$:
$G=A*_P P*_{P'} B' *_{Q'} Q*_{Q} C$. Although $\C'$ has two more vertices
than $\C$ ,
$\C'$ is not reduced. And if we collapse edges of $\C'$ to obtain a reduced decomposition
we get
$G=A*_{P'} B'*_{Q'} C$, which has the same number of vertices  as $\C$.
\end{example}

As this example shows Theorem \ref{thm,bestvinafeighn} can not
control a sequence of reduced graph decomposition which is
obtained by elementary unfoldings. But we have another
accessibility result to control this, which we prove later.

As we said, we now produce a sequence of refinements $\Delta_n$
using the splittings of $G$ along $C_n$. We may assume that the
first splitting is the splitting along $C$, with which we already
constructed $\Delta_1$. We now refine $\Delta_1$ using the
splitting along $C_2$. Same as $\Lambda$, all edge groups and all
enclosing vertex groups of $\Delta_1$ are elliptic on $T_2$, the
Bass-Serre tree of the splitting along $C_2$. As before, we let a
non-enclosing vertex group of $\Delta_1$, $U$, act on $T_2$ and
obtain a graph decomposition of $U$, which we substitute for the
vertex labelled by $U$ in $\Delta_1$. We do this for all
non-enclosing vertex groups of $\Delta_1$, then we apply Prop
\ref{prop,modification}. If the resulting graph decomposition is
an elementary unfolding of $\Delta_1$, then we put
$\Delta_2=\Delta_1$. Otherwise, if the graph decomposition is not
reduced, then we collapse one edge at a vertex where it is not
reduced, and obtain a reduced graph decomposition of $G$, which we
denote $\Delta_2$. $\Delta_2$ has following properties:
\begin{enumerate}
\item
$\Delta_2$ is a refinement of $\Delta_1$.
$\Delta_2$ is identical to $\Delta_1$,  or has more
vertices.
\item
Each edge of $\Delta_2$ gives a minimal
splitting of $G$ along a slender group.
The edge group  is a subgroup
of a conjugate of either an edge group of $\Delta_1$ or $C_2$.
\item
Each maximal enclosing group is a conjugate of some
vertex group of $\Delta_2$.
\item
After, if necessary, performing an elementary unfolding to
$\Delta_2$,
each vertex group is elliptic on $T_2$, and also $T_1$, the
Bass-Serre tree of the splitting along $C_1$.
\end{enumerate}

We repeat the same process; refine $\Delta_n$ using the splitting
along $C_n$ to $\Delta_{n+1}$. Because of the property 1 in the
above list, by Theorem \ref{thm,bestvinafeighn}, there exists a
number, $N$, such that if $n \ge N$, then $\Delta_{n+1}$ is equal
to or an elementary unfolding of $\Delta_n$.
 Let's denote $\Delta_N$ by $\Delta$.
We state some properties of $\Delta$.

\begin{prop}[JSJ decomposition for minimal splittings with elementary unfoldings]
\label{prop,jsjwithunfolding}
Let $G$ be a finitely presented group.
Let $G=A_n*_{C_n} B_n$ (or $A_n*_{C_n}$) be all minimal splittings
of $G$ along slender groups, $C_n$, which are elliptic-elliptic
with respect to any minimal splittings of $G$ along slender subgroups.
Let $T_n$ be their Bass-Serre trees.
Then there exists a graph decomposition, $\Delta$, of $G$ such
that
\begin{enumerate}
\item[1,2,3.] Same as the properties 1,2,3 of 
Prop \ref{prop,jsj.hyp-hyp.minimal}.
\item[4.] For each $n$, there exists an elementary unfolding of
$\Delta$ such that each vertex group is elliptic on  $T_n$.
\item[5.] Let $G=A*_C B$ or $G=A*_C$ be a minimal splitting along
a slender group $C$ which is hyperbolic-hyperbolic with respect to
some minimal splitting along a slender group, and $T_C$ its
Bass-Serre tree. Then,
\begin{enumerate}
\item same as \ref{item,hyperbolic.enclosing} of Prop
\ref{prop,jsj.hyp-hyp.minimal}.
\item
There exists an elementary unfolding
of $\Delta$, at non-enclosing vertex groups,
 such that  all vertex groups in the elementary unfolding
except for $S$ are elliptic on $T_C$, the Bass-Serre tree of the
splitting along $C$.
\end{enumerate}
\end{enumerate}
\end{prop}

\begin{rem}
In fact we do not need an elementary unfolding in the properties 4
and 5 in the proposition, if we construct a more refined $\Delta$.
We show this in Theorem \ref{thm,jsj.minimal}.
\end{rem}

\proof We already know 1,2,3,4 from the way we constructed
$\Delta$. Also the property 5(a) is immediate from Prop
\ref{prop,jsj.hyp-hyp.minimal}. To show 5(b), let $U$ be a
non-enclosing vertex group of $\Delta$ which is not elliptic on
$T_C$. If such vertex does not exist, we are done. As usual,
letting $U$ act on $T_C$, we obtain a graph decomposition,
${\mathcal U}$, of $U$, then we substitute this for $U$ in
$\Delta$ to obtain a refinement, $\Delta'$, of $\Delta$ whose
edges give minimal splittings, after we apply Prop
\ref{prop,modification} if necessary. But this resulting
decomposition has to be an elementary unfolding of
$\Delta=\Delta_N$, because otherwise we must have refined
$\Delta_N$ further when we constructed $\Delta$.
\qed

\subsection{Elementary unfolding and accessibility}
As we said in the remark after Prop \ref{prop,jsjwithunfolding},
we do not need elementary unfoldings.
But as we saw in Example \ref{example,unfolding}, Theorem \ref{thm,bestvinafeighn}
can not control a sequence of elementary unfoldings because they are
not reduced.
We prove another accessibility result.
This result was suggested to us by Bestvina.
The argument is similar to the one used by Swarup in a proof of
Dunwoody's accessibility result (\cite {Sw}).

\begin{prop}[Intersection accessibility]
\label{prop,intersection}
Let $G$ be a finitely presented group. Suppose
$\C_i$ is a sequence of graph decompositions of $G$
such that all edge groups are slender.
Suppose for any $i$, $\C_{i+1}$ is obtained from
 $\C_i$ by an elementary unfolding. Then there is a graph decomposition $\C$
of $G$ with all edge groups slender such that
for any $i$, $\C $ is a refinement of $\C_i$.
\end{prop}

\proof
We define a partial order on the set of graph of groups
decompositions of $G$.
We say that $\Gamma < \Lambda$
if all vertex groups of $\Gamma$ act elliptically
on the Bass-Serre tree corresponding to $\Lambda$,
in other words, $\Gamma$ is a refinement of $\Lambda$.
We have $\C_{i+1}< \C_i$ for all $i$.

We can apply Dunwoody's tracks technique to obtain a graph of groups
decomposition $\Gamma $ such that  $\Gamma <\Gamma
_i$ for
all $i$. We describe briefly how this is done: Let $K$ be a presentation
complex for $G$. Without loss of generality we assume that $K$
corresponds to
a triangular presentation.
Let $T_i$ be the Bass-Serre tree of $\C _i$. As we noted earlier there are maps $\phi _i:
T_{i+1}\to T_i$ obtained by collapsing some edges.
We
choose a sequence of  points $(x_i)$ such that $x_i$ is a midpoint of an
edge
and $\phi (x_{i+1})=x_i$.

We will define maps $\alpha _i:K\to T_i$. Each oriented edge of $K$
corresponds
to a generator of $G$. Given an edge $e$ corresponding to an element
$g\in G$
we map it by a linear map to the geodesic joining $x_i$ to $gx_i$. We
extend
linearly this map to the 2-skeleton of $K$. A track is a preimage  of
a vertex of $T_i$ under this map.
We note that the tracks we obtain from $T_i$ are a subset of the tracks
obtained from $T_{i+1}$ (or to be more formal each track obtained from
$T_i$ is
`parallel' to a track obtained from $T_{i+1}$). We remark that
for each $i$ the tracks obtained from $\alpha _i$ give rise to a decomposition
$\Gamma _i'$ of $G$. $\C _i$ is obtained from $\C _i'$ by subdivisions
and foldings.

Since $G$ is finitely presented, so that
$K$ is compact, there is a $\lambda (G)$ such that there are at most $\lambda (G)$
non-parallel tracks we conclude that there is an $n$ such that each
track
obtained from $T_k$ ($k>n$) is parallel to a track obtained from $T_n$.
We can then take as $\Gamma_{\infty} $ the graph
of groups decomposition corresponding to the tracks obtained from
 $T_n$. It follows that $\Gamma_i > \C_{\infty}$.
Put $\Gamma=\Gamma_{\infty}$.
\qed

\subsection{JSJ-decomposition along slender groups}
We state one of our main theorems.
\begin{thm}[JSJ-decomposition for minimal splittings
along slender groups]
\label{thm,jsj.minimal}

Let $G$ be a finitely presented group.
Then there exists a graph decomposition, $\Gamma$, of $G$
such that
\begin{enumerate}
\item[1,2,3.] same as the properties 1,2,3 of Prop \ref{prop,jsj.hyp-hyp.minimal}.
\item[4.]
Let $G=A*_C B$ or $A*_C$ be a minimal splitting along
a slender group $C$, and $T_C$ its Bass-Serre tree.
\begin{enumerate}
\item
If it is elliptic-elliptic with respect to all
minimal splittings of $G$ along slender groups, then
all vertex groups of $\Gamma$ are elliptic on $T_C$.
\item
If it is hyperbolic-hyperbolic
with respect to some minimal splitting of $G$ along a slender
group, then there is an enclosing vertex group, $S$, of $\Gamma$
which contains a conjugate of $C$ and
the property  \ref{item,hyperbolic.enclosing} of Prop
\ref{prop,jsj.hyp-hyp.minimal} holds for $S$.
All vertex groups except for $S$ of $\Gamma$ are elliptic on $T_C$.

In particular, there is a graph decomposition, ${\mathcal S}$, of $S$
whose edge groups are in conjugates of $C$, which we can
substitute for $S$ in $\Gamma$ such that all vertex groups
of the resulting refinement of $\Gamma$ are elliptic on $T_C$.
\end{enumerate}
\end{enumerate}
\end{thm}

\proof Let $\Delta$ be the graph decomposition of $G$ which we
have constructed for Prop \ref{prop,jsjwithunfolding}. We will
obtain $\Gamma$ as a refinement of $\Delta$ at non-enclosing
vertex groups. Let $G=A_n*_{C_n} B_n$ (or $A_n*_{C_n}$) be all
minimal splittings of $G$ along slender groups, $C_n$, which are
elliptic-elliptic with respect to any minimal splitting of $G$
along slender subgroups, and $T_n$ their Bass-Serre trees. We have
defined a process to refine a graph decomposition using this
collection to obtain a sequence $\{\Delta_n\}$ for Prop
\ref{prop,jsjwithunfolding}. We apply nearly the same process to
$\Delta$ again using the splittings along $C_n$, and produce a
sequence $\{\Gamma_n\}$. The only difference is that we do not
make a graph decomposition reduced in each step. Let $T_n$ be the
Bass-Serre tree of the splitting along $C_n$. To start with, put
$\Gamma_0=\Delta$. Letting all vertex groups act on $T_1$, we
obtain graph decompositions, then substitute them for the
corresponding vertex groups in $\Gamma_0$, which is $\Gamma_1$.
$\Gamma_1$ is an elementary unfolding of $\Gamma_0$, because
otherwise, we must have refined $\Delta$ farther in the proof of
Prop \ref{prop,jsjwithunfolding}. Note that $\Gamma_1$ is not
reduced, but we do not collapse any edges. We repeat the same
process; we let all vertex groups of $\Gamma_1$ act on $T_2$,
substitute those graph decompositions for the corresponding vertex
groups in $\Gamma_1$. The resulting non-reduced graph
decomposition is $\Gamma_2$, and so on. In this way, we obtain a
sequence of graph decompositions $\Gamma_n$ such that
$\Gamma_{n+1}$ is an elementary unfolding of $\Gamma_n$. We remark
that in each step enclosing vertex groups stay unchanged since
they are elliptic on all $T_n$. Note that each $\Gamma_n$
satisfies the properties 1,2,3 and 5(a)i of Prop
\ref{prop,jsj.hyp-hyp.minimal}.

Suppose that  there exists $N$ such that for any $n \ge N$,
$\Gamma_n=\Gamma_{n+1}$.
Then $\Gamma_N$ satisfies the properties 4 and 5(b) of
Prop \ref{prop,jsjwithunfolding} as well without
elementary unfoldings, so that
the property 4 of the theorem follows.
Putting $\Gamma=\Gamma_N$, we obtain a desired $\Gamma$.

If such $N$ does not exist, then we apply Prop
\ref{prop,intersection} to our sequence and obtain a graph
decomposition which is smaller (or equal to), for the order
defined in Prop \ref{prop,intersection}, than all $\Gamma_n$.
Let's take a minimal element, $\Gamma$, with respect to our order.
Such a decomposition exists by Zorn's lemma. $\Gamma$ is the
decomposition that we look for, because if we apply the process to
refine $\Gamma$ using the sequence of decompositions along $C_n$
as before, nothing happens, because $\Gamma$ is minimal in our
order. It follows that $\Gamma$ satisfies all the properties.
\qed

We call a graph decomposition of $G$ we obtain in Theorem
\ref{thm,jsj.minimal} a {\it JSJ decomposition of $G$ for
splittings along slender groups}. We will prove that $\Gamma$ has
the properties stated in this theorem not only for minimal
splittings of $G$ along slender groups, but also non-minimal
splittings as well in Theorem \ref{thm,jsj}.

\begin{cor}[Uniqueness of JSJ decomposition]\label{cor,uniqueness.jsj.minimal}
Let $G$ be a finitely presented group.
Suppose a graph decomposition, $\Gamma$, of $G$
satisfies the properties 2 and 4(a) of Theorem
\ref{thm,jsj.minimal}.
\begin{enumerate}
\item Suppose  $\Gamma'$ is a graph decomposition of $G$ which
satisfies the properties 2 and 4(a) of Theorem
\ref{thm,jsj.minimal}. Then all vertex groups of $\Gamma'$ are
elliptic on the Bass-Serre tree for $\Gamma$. \item $\Gamma$
satisfies the property 3 of Theorem \ref{thm,jsj.minimal}. \item
$\Gamma$ satisfies the property 4(b) if $G$ does not split along
an infinite index subgroup of $C$.
\end{enumerate}
\end{cor}

\proof 1. Let $T$ be the Bass-Serre tree for $\Gamma$. Let $V$ be
a vertex group of $\Gamma'$. Let $e$ be an edge of $\Gamma$ with
edge group $E$, and $T_e$ the Bass-Serre tree of the splitting of
$G$ along $E$ which the edge $e$ gives. To show $V$ is elliptic on
$T$, it suffices to show that it is elliptic on $T_e$ for all $e$.
This splitting along the slender group $E$ is minimal, and
elliptic-elliptic with respect to any minimal splitting of $G$
along a slender group by the property 2 of $\Gamma$. Therefore, by
the property 4(a) of $\Gamma'$, all vertex groups of $\Gamma'$ are
elliptic on $T_e$. In particular $V$ is elliptic on $T_e$, so it
is elliptic on $T$.

2. Let $T$ be the Bass-Serre tree of $\Gamma$. Let $S$ be the
maximal enclosing vertex group in a maximal enclosing
decomposition, $\Lambda$, of $G$. We first show that $S$ is
elliptic on  $T$. To show it, as usual, we use Cor
\ref{cor,elliptic}. Let $\Sigma$ be the $2$-orbifold for $S$, and
$s$ an essential simple closed curve or a segment on $\Sigma$.
Cutting $\Sigma$ along $s$, we obtain a splitting of not only $S$
but also $G$ along the slender group, $C$, represented by $s$.
This splitting is minimal by Prop \ref{prop,dual-minimal}. By Cor
\ref{cor,elliptic}, it suffices to show that the group $C$ is
elliptic on $T$. By the property 1 of $\Gamma$, the pair of the
splittings of $G$ along $E$ (from the previous paragraph) and $C$
is elliptic-elliptic. Therefore, the group $C$ is elliptic on
$T_E$, so that it is elliptic on $T$ as well since the edge $e$
was arbitrary.

We already know that $S$ is in a conjugate of a vertex group, $V$,
of $\Gamma$. We want to show that indeed $S$ is a conjugate of
$V$. Let $T_{\Lambda}$ be the Bass-Serre tree of the maximal
enclosing decomposition $\Lambda$  with $S$. It suffices to show
that $V$ is elliptic on $T_{\Lambda}$ to conclude that $S$ is a
conjugate of $V$, because $S$ is the only vertex group of
$\Lambda$ which can contain a conjugate of $V$. This is because
all edge groups adjacent to $S$ are peripheral subgroups, so they
are proper subgroups of $S$. Let $d$ be an edge of $\Lambda$ with
the edge group $D$. The splitting of $G$ along $D$ which the edge
$d$ gives
 is minimal since
$\Lambda$ is an enclosing decomposition, and elliptic-elliptic
with respect to any minimal splitting along a slender group
since $\Lambda$ is maximal.
By the property 3 (a) of $\Gamma$, all vertex groups
of $\Gamma$ are elliptic on $T_D$, the Bass-Serre tree
of the splitting along $D$. Since the edge $d$ was arbitrary,
all vertex groups of $\Gamma$ are elliptic
on $T_{\Lambda}$, in particular, so is $V$.

3. Let $I$ be a maximal set of hyperbolic-hyperbolic minimal
splittings of $G$ along slender groups which contains the
splitting along $C$. Let $\Lambda$ be a maximal enclosing
decomposition of $G$ for $I$ with enclosing vertex group $S$. Let
$\Sigma$ be the $2$-orbifold for $S$. We can assume that there is
an essential simple closed curve or a segment, $s$, on $\Sigma$
such that by cutting $\Sigma$ along $s$ we obtain a splitting of
$G$ along $C$. This is because when we construct $\Lambda$ for $I$
using a sequence of graph decomposition, we can start with the
splitting along $C$. Although we do not know in general if this
splitting is the same as the one we are given, it is the case
under our extra assumption. By the property 3 of $\Gamma$, a
conjugate of $S$ is a vertex group of $\Gamma$, which therefore
contains a conjugate of $C$. No other vertex group of $\Gamma$
contains a conjugate of $C$ because if it did, then an edge group
of $\Gamma$ has to contain a conjugate of $C$, which is a
contradiction since the splitting along $C$ is
hyperbolic-hyperbolic while $\Gamma$ has property 2.

Let $w$ be the vertex of $\Gamma$ with the vertex group, $W$,
which is a conjugate of $S$. Let $v$ be a vertex of $\Gamma$ with
 vertex group, $V$, such that $v \not=w$. We want to show that
$V$ is elliptic on $T_C$, the Bass-Serre tree of the splitting
$G=A*_C B,$ or $ A*_C$, which we know is obtained by cutting
$\Sigma$ along $s$. We first claim that $V$ is elliptic on
$T_{\Lambda}$. This is because each edge of $\Lambda$ gives a
minimal splitting which is elliptic-elliptic since $\Lambda$ is
maximal, so that $V$ is elliptic on $T_{\Lambda}$ by the property
4(a). (Use it to each edge decomposition of $\Lambda$). Therefore
$V$ is in a conjugate of some vertex group, $U$, of $\Lambda$. If
$U$ is not $S$, we are done, because then $U$ is elliptic on
$T_C$. We have used that the original splitting along $C$ is
identical to the one we obtain by cutting along $s$. Suppose
$U=S$, then $V$ is in a conjugate of $W$. By Bass-Serre theory,
this means that there is an edge in $\Gamma$ adjacent to $v$ whose
edge group is $V$. It then follows from the property 2 for
$\Gamma$ that the edge group $V$ is elliptic on $T_{\Lambda}$. The
proof is complete. \qed

As we said, $\Gamma$ indeed can deal with non-minimal splittings
of $G$ along slender groups as well.

\begin{thm}[JSJ decomposition for splittings along slender groups]
\label{thm,jsj} Let $G$ be a finitely presented group, and let
$\Gamma$ be the graph decomposition we obtain in Theorem
\ref{thm,jsj.minimal}. Let $G=A*_C B, A*_C$ be a splitting along a
slender group $C$, and $T_C$ its Bass-Serre tree.
\begin{enumerate}
\item
If the group $C$ is elliptic with respect to any minimal splitting of $G$ along
a slender group, then all vertex groups of $\Gamma$ are
elliptic on $T_C$.
\item
Suppose the group $C$ is hyperbolic with respect to
some minimal splitting of $G$ along
a slender group. Then
\begin{enumerate}
\item all non-enclosing vertex groups of $\Gamma$ are elliptic on
$T_C$. \item For each enclosing vertex group, $V$, of $\Gamma$,
there is a graph decomposition of $V$, ${\mathcal V}$, whose edge
groups are in conjugates of $C$, which we can substitute for $V$
in $\Gamma$ such that if we substitute for all enclosing vertex
groups of $\Gamma$ then all vertex groups of the resulting
refinement of $\Gamma$ are elliptic on $T_C$.
\end{enumerate}
\end{enumerate}
\end{thm}

\proof 1. If the splitting along $C$ is minimal, then nothing to
prove (Theorem \ref{thm,jsj.minimal}). Suppose not. Apply Prop
\ref{prop,modification} to the splitting and obtain a refinement,
$\Lambda$, such that each edge, $e$, of $\Lambda$ gives a minimal
splitting of $G$ along a slender group, $E$, which is a subgroup
of a conjugate of $C$. Let $G=P*_E Q$ (or $ P*_E$) be the
splitting along $E$ which the edge $e$ gives. Let $T_E$ be its
Bass-Serre tree. Let $T_{\Lambda}$ be the Bass-Serre tree of
$\Lambda$. We want to prove that each vertex group, $V$, of
$\Gamma$ is elliptic on $T_{\Lambda}$, which implies that $V$ is
elliptic on $T_C$, since $\Lambda$ is a refinement of the
splitting along $C$. By Bass-Serre theory, it suffices to prove
that $V$ is elliptic on $T_E$. By our assumption, the group $C$ is
elliptic with respect to any minimal splitting of $G$ along a
slender group, so that so is $E$ since it is a subgroup of a
conjugate of $C$. Therefore, the minimal splitting $G=P*_E Q$ (or
$P*_E$) is elliptic-elliptic with respect to any minimal splitting
of $G$ along a slender group, so that by the property 4(a),
Theorem \ref{thm,jsj.minimal} of $\Gamma$, all vertex groups of
$\Gamma$, in particular $V$, are elliptic on $T_E$.

2. If the given splitting along $C$ is minimal, then nothing to
prove because we have 4(b), Theorem \ref{thm,jsj.minimal}. Suppose
it is not minimal, and apply Prop \ref{prop,modification} to
obtain a refinement, $\Lambda$ such that each edge gives a minimal
splitting of $G$. All edge groups of $\Lambda$ are in conjugates
of $C$.  Let $T_{\Lambda}$ be its Bass-Serre tree. Then all
non-enclosing vertex groups of $\Gamma$ are elliptic on
$T_{\Lambda}$. The argument is similar to the case 1 above and the
proof for 3(b), Theorem \ref{thm,jsj.minimal}. We omit details.
Let $V$ be an enclosing vertex group of $\Gamma$. Letting $V$ act
on $T_{\Lambda}$, we obtain a graph decomposition, ${\mathcal V}$,
of $V$ such that all edge groups are in conjugates of $C$. Because
each edge of $\Lambda$ gives a minimal splitting of $G$ along a
slender group, by the property 1, Theorem \ref{thm,jsj.minimal}
for $\Gamma$, all edge groups of $\Gamma$ are elliptic on
$T_{\Lambda}$. Therefore we can substitute ${\mathcal V}$ for $V$
in $\Gamma$. If we substitute for all enclosing vertex groups in
this way, we obtain the desired refinement of $\Gamma$. \qed

One can interpret our theorems using the language of foldings.
What we show is that if $\Gamma $ is the JSJ-decomposition of a
finitely presented  group and $A*_CB$ (or $A*_C$) is a splitting of
$G$ over a slender group $C$ with Bass-Serre tree $T_C$ then we
can obtain a graph decomposition $\Gamma '$ from $\Gamma $ such
that all vertex groups of $\Gamma '$ act elliptically on $T_C$.
Let's call $T$ the Bass-Serre tree of $\Gamma '$. Since all vertex
groups of $\Gamma '$ fix vertices of $T_C$ we can define a
$G$-equivariant simplicial map $f$ from a subdivision of $T$ to
$T_C$. To see this pick a tree $S\subset T$ such that the
projection from $S$ to $\Gamma $ is bijective on vertices. If
$v\in S^0$ pick a vertex $u\in T_C$ such that $Stab(v)\subset
Stab(u)$. Define then $f(v)=u$. We extend this to the edges of $S$
by sending the edge joining two vertices to the geodesic joining
their images in $T_C$. This can be made simplicial by subdividing
the edge. Finally extend this map equivariantly on $T$. It follows
that the splitting over $C$ can be obtained from $\Gamma $ by
first passing to $\Gamma '$ and then performing a finite sequence
of subdivisions and foldings. In this sense $\Gamma $ 'encodes'
all slender splittings of $G$.

\section{Final remarks}
 For producing the
JSJ-decomposition we did not put any restriction on $G$; in
particular we did not assume that $G$ does not split over groups
`smaller' than the class considered. One of the difficulties in
this is that there is no natural `order' on the set of slender
groups. Otherwise one could work inductively starting from the
`smallest' ones. This is the reason we introduced the notion of
minimal splittings.

We remark that the situation is simpler if one restricts one's
attention to polycyclic groups as one can 'order' them.

It is a natural question whether there is a JSJ-decomposition over
small groups. Our results (in particular proposition
\ref{prop,enclosinggroup}) might prove useful in this direction.
The main difficulty for generalizing it to an arbitrary number of
small splittings to produce a JSJ-decomposition over small groups
is that the edge groups of the decomposition in prop.
\ref{prop,enclosinggroup} are not small in general. So one can not
apply induction in the case of small splittings.

  We note however that if a
JSJ-decomposition over small groups exist its edge groups are not
small. We illustrate this by the following example.

\begin{example}
Let's denote by $A$ the group given by
$A=\langle a,t,s|tat^{-1}=a^2,sas^{-1}=a^2 \rangle$ and let $H$ be an
unsplittable group containing $F_2$, e.g. $SL_3(\mathbb{Z})$.

Let's consider a complex of groups $G(X)$ with underlying complex
a sphere obtained by gluing two squares along their boundary. We
label all 4 vertices (0-cells) by $A\times H$. We label 2-cells
and 1-cells by an infinite cyclic group $\langle c \rangle$. In one of the
squares all maps from the group of a 2-cell to a group of a 1-cell
are isomorphisms while all maps to the group of 0-cell send $c$ to
$a$.

We describe now the maps in the second square $\tau $. Let
$e_1,e_2$ be two adjacent edges (1-cells) of the square. Let
$a_{12}$ be the common vertex of $e_1,e_2$ and let $a_1$ be the
other vertex of $e_1$ and $a_2$ the other vertex of $e_2$. Let
finally $b$ be the fourth vertex of the square.

The monomorphisms $\psi _1: G_{\tau }\to G_{e_1}, \psi _2: G_{\tau
}\to G_{e_2} $ are given by $c\to c^2$. All other maps are defined
as in the first square. To satisfy condition 3 of the definition
of a complex of groups (see subsec. 4.1) we define the 'twisting'
element $g_{e,f}$ for two composable edges $e,f$ as follows:

For each vertex of the square there are two pairs of composable
edges from the barycenter of $\tau $ to the vertex.
For the vertex
$b$ we put $g_{e,f}=1$ for the first pair and $g_{e,f}=t^{-1}s$
for the second pair. We remark that $t^{-1}s$ commutes with $a$ so
condition 3 is satisfied.
For the vertex $a_1$ for the pair of composable edges that
corresponds to an isomorphism from $G_{\tau }$ to $G_{a_1}$ we put
$g_{e,f}=1$ and for the other pair we put $g_{e,f}=t$.
Similarly for the vertex $a_2$ for the pair of composable edges
that corresponds to an isomorphism from $G_{\tau }$ to $G_{a_2}$
we put $g_{e,f}=1$ and for the other pair we put $g_{e,f}=t$.
Finally for the vertex $a_{12}$ for one pair we put $g_{e,f}=t$
and for the other $g_{e,f}=s$.

 It is now a straightforward
computation to see that this complex is developable. In fact to
show this here it is enough to show that links of vertices do not
contain simple closed curves of length 2. Let $v$ be a vertex of
$X$. The link of $v$, $Lk\,v$ has as set of vertices the pairs
$(g\psi _a(G_{i(e)}),e)$ where $g\in G_v$ and $e\in E(X)$ with
$t(e)=v$. The set of edges of the barycentric subdivision of
$Lk\,v$ is the set of triples $(g\psi _{ef} (G_{i(f)}),e,f)$ where
$e,f$ are composable edges in $E(X)$ with $t(e)=v$. The initial
and terminal vertices of an edge are given by:
$$i(g\psi _{ef} (G_{i(f)}),e,f)=(g\psi _{ef} (G_{i(f)}),ef) $$
$$i(g\psi _{ef} (G_{i(f)}),e,f)=(gg_{e,f}^{-1}\psi _{e} (G_{i(e)}),e) $$
Our choice of twisting elements now insures that there are no
curves of length 2 in the link. To see this notice that if e.g.
one assigns $g_{e,f}=1$ to all pairs of composable edges in $E(X)$
with $t(e)=b$ then condition 3 of the definition of the complex of
groups is satisfied but now the link has a simple closed curve of
length 2. Our choice of non-trivial twisting element insures that
there are no length 2 curves in the link.

Let's denote by $G$ the fundamental group of $G(X)$.
We remark that the two simple closed curves perpendicular at the
midpoints of $e_1,e_2$ give rise to two small splittings of $G$
over $BS(1,2)=\langle x,y|xyx^{-1}=y^2 \rangle $. Note also that this pair of
splittings is hyperbolic-hyperbolic.

 We claim that this complex gives the
JSJ-decomposition of $G$ over small groups.
Let $\tilde X$ be a complex on which $G$ acts with quotient
complex of groups $G(X)$. Let $T$ be the Bass-Serre tree of a
 splitting of $G$ over a small group $C$. Then all vertex stabilizers of $\tilde X$
fix vertices of $T$. This implies that there is a $G$ -equivariant
map $f:\tilde X \to T$. The preimage of a midpoint of an edge of
$T$ is a graph (a tree) in $\tilde X $ projecting to an essential
simple closed curve on $X$ corresponding to a splitting of $G$
over a conjugate of $C$. In other words this complex gives us a
JSJ-decomposition for $G$. Now one of the edge groups of this
decomposition has the presentation $\langle s,a|sa^2s^{-1}=a^2 \rangle$ (it is
the edge corresponding to the vertex lying in both $e_1,e_2$),
which clearly is not a small group. We remark that 2 edge groups
are labelled by $BS(1,2)$ and one by $\mathbb{Z}\times
\mathbb{Z}$, so they are small.
\end{example}

\end{document}